\documentclass{amsart}                     

\usepackage{framed}
\usepackage{amssymb}
\usepackage{amsmath}
\usepackage{url}
\usepackage{xspace}
\usepackage{graphics,graphicx}

\usepackage{multirow}
\usepackage{alltt}
\usepackage{a4wide}
\usepackage{latexsym}
\usepackage{epsfig}
\usepackage{amssymb}
\usepackage{color}
\usepackage{textcomp} 
\let\proof\relax

\usepackage{amsthm}

\newtheorem{definition}{Definition}
\newtheorem{theorem}{Theorem}
\newtheorem{lemma}{Lemma}
\newtheorem{remark}{Remark}
\newtheorem{proposition}{Proposition}

\newcommand{\ub}{\,\mbox{$\bullet$}\,}
\newcommand{\ob}{\,\mbox{\bf\texttt{[}}\,}
\newcommand{\cb}{\,\mbox{\bf\texttt{]}}\,}
\newcommand{\op}{\,\mbox{\bf\texttt{(}}\,}
\newcommand{\cp}{\,\mbox{\bf\texttt{)}}\,} 
\def\cA{\mathcal{A}}
\newcommand{\cF}{\mathcal{F}} 
\newcommand{\cG}{\mathcal{G}}
\def\NN{\mathbb{N}}

\newcommand{\cacher}[1]{}

\newcommand{\mfold}{\mbox{\tt mfold}\xspace}
\newcommand{\rnafold}{\mbox{\tt RNAfold}\xspace}

\chardef\other=12

\def\mmakeactive#1{\catcode`#1=\active\ignorespaces}

{
\mmakeactive\
\gdef\obeywhitespace{%
  \mmakeactive\^^M %
  \let^^M=\NewLine %
  \aftergroup\removebox %
  \obeyspaces %
}}

\def\NewLine{\par\indent}
\def\removebox{\setbox0=\lastbox}

\def\|{|}

\title{Combinatorics of locally optimal RNA secondary structures}
\author[\'E. Fusy and P. Clote]{\'Eric Fusy$^{*}$ \and Peter Clote$^{\dagger}$}
\thanks{$^{*}$ LIX, \'Ecole Polytechnique, Palaiseau, France. fusy@lix.polytechnique.fr.\\ 
$^{\dagger}$ Department of Biology, Boston  College, Chestnut Hill, MA 02467, 
USA. clote@bc.edu.}

\begin{document}
\maketitle

\begin{abstract}
It is a classical result of Stein and Waterman
that the asymptotic number of RNA secondary structures  is
$1.104366 \cdot n^{-3/2} \cdot 2.618034^n$. Motivated by 
the kinetics of RNA secondary structure formation,
we are interested in determining the asymptotic number of 
secondary structures that are {\em locally optimal}, with respect
to a particular energy model. In the Nussinov energy model, where
each base pair contributes $-1$ towards the energy of the structure,
locally optimal structures are exactly the {\em saturated} structures,
for which we have previously shown that asymptotically, there are
$1.07427\cdot n^{-3/2} \cdot 2.35467^n$ many saturated structures for
a sequence of length $n$. In this paper, we consider the {\em base
stacking energy model}, a mild variant of the Nussinov model, where
each stacked base pair contributes $-1$ toward the energy of the structure.
Locally optimal structures with respect to the base stacking energy model
are exactly those secondary structures, whose stems cannot be extended.
Such structures were first considered by Evers and Giegerich, who
described a dynamic programming algorithm to
enumerate all locally optimal structures.
In this paper, we apply methods from enumerative combinatorics 
to compute the asymptotic number of such structures. Additionally,
we consider analogous combinatorial problems for secondary structures with
annotated single-stranded, stacking nucleotides (dangles).
\end{abstract}

\section{Introduction}

Historically, the development of combinatorics 
for RNA secondary structures \cite{steinWaterman,waterman:SecStr}
has been intimately related to the development of
{\em algorithms} for RNA minimum free energy (MFE) secondary
structure \cite{zukerStiegler,zuker:mfoldWebserver,hofacker:RNAfoldServer}.
In particular,  {\em counting} the number of secondary
structures for sequence of length $n$ is essentially equivalent to
computing the Boltzmann partition function, defined
by $Z= \sum_{S} \exp(-E(S)/RT)$, where the sum
is taken over all secondary structures $S$, the energy of $S$ is
denoted by $E(S)$, $R \approx 1.959$ cal/mol 
is the universal gas constant, and $T$
absolute temperature.\footnote{If the energy $E(S)=0$ or if the
temperature $T=+\infty$, then the partition function is exactly
equal to the number of secondary structures.}

Complex analysis is used to obtain the asymptotic enumeration
results described in this article and related articles mentioned in
the introduction.
In particular, given a complex generating function 
$f(z) = \sum a_n z^n$,
it is well-known from introductory complex analysis that $f$
converges in a circular region about the point of expansion out to the 
dominant, or nearest, singularity $r$, and thus the asymptotic order
of magnitude of $a_n$ is approximately $r^{-n}$.  Darboux's 
theorem\footnote{Jean Gaston Darboux (1842--1917).}
\cite{siamreviewcombinatorics,henrici:Combinatorics} states that
if $f(z) = \sum_{n=0}^{\infty} (r-z)^{\alpha} L(z)$, where $r>0$,
$\alpha$ is not a positive integer, and $L$ is analytic in a disk of
radius greater than $r$, then $\alpha_n \sim r^{\alpha-n} n^{\alpha-1} 
\frac{L(r)}{\Gamma(-\alpha)}$. This result was generalized by
Bender \cite{bendertheoremoriginal}, corrected by
Meir and Moon \cite{meirMoon}, and further extended by 
Flajolet and Odlyzko \cite{flajolet-odlyzko} and by Drmota,
Lalley and Woods (each of the latter worked independently) -- see the
exposition in \cite{FlSebook} for discussion and references. 

In \cite{steinWaterman}, Stein and Waterman proved that
the asymptotic number of secondary structures is
$1.104366 \cdot n^{-3/2} \cdot 2.618034^n$.  Since that time, a number
of additional results on the combinatorics of RNA structures have
been obtained.
In \cite{hofacker99a}, Hofacker et al. derived a number of asymptotic
results on the number of structures, expected number of base pairs, etc.
for RNA secondary structures. Observing a correspondence with involutions,
Haslinger and Stadler \cite{Haslinger-Stadler-1999} provided an upper bound
on the number of {\em bi-secondary} structures, i.e. structures having
non-nested pseudoknots that can be presented as a union
$S=S_1 \cup S_2$ of disjoint secondary structures, and Rodland
\cite{rodland:pseudoknots} studied the asymptotic number of
a number of classes of pseudoknotted structures. Building on a remarkable
and pioneering paper of Harer and Zagier \cite{harerZagier}, 
Vernizzi et al.  \cite{Vernizzi.prl05} classified
pseudoknotted RNA structures according to topological genus $g$, and
then applied the work of Harer and Zagier to obtain recurrence relations for 
the number of pseudoknotted structures of genus $g$.  In 
\cite{Saule.jcb11}, Saule et al. provided a summary table of the asymptotic
number of pseudoknotted structures structures, with respect to various
allowed pseudoknots, and established the asymptotic number
of pseudoknotted structures, with no restriction. In \cite{Li.jmb12},
Li and Reidys determined the asymptotic number of hybridizations of
two interacting RNA structures. Moving away from counting the number
of structures, Yoffe et al. \cite{Yoffe.nar11}
and Clote et al. \cite{Clote.jmb11} determined the 
asymptotic expected distance between the $5'$ and $3'$ ends of RNA sequence, 
where the $5'$ to  $3'$ distance of a given structure $S$ on sequence
$s_1,\ldots,s_n$ is defined as the 
minimum number of backbone or base-pairing edges in a minimum length
path from $s_1$ to $s_n$.

In \cite{Clote.jcb06}, Clote computed the asymptotic number 
$1.07427\cdot n^{-3/2} \cdot 2.35467^n$ of
{\em saturated} structures, defined by Zuker \cite{zuker:1986}
as those for which no base pair
can be added without violating the definition of secondary structure.
In \cite{Clote.jbcb09}, Clote et al. provided another proof for the 
asymptotic number of saturated structures, which additionally yielded the
asymptotic expected number of base pairs $0.337361 \cdot n$ for
saturated structures.
An overview of methods for RNA enumerative combinatorics is given in
Lorenz et al. \cite{Lorenz.jcb08}, where additionally it is shown that
the asymptotic number of {\em shapes} of secondary structures for
a length $n$ sequence is $2.44251 \cdot n^{-3/2} \cdot 1.32218^{n}$.\footnote{The {\em shape} of a secondary structure was defined by Voss et al.
\cite{giegerich:shapesProbAnal} to represent its branching topology; for
instance, the shape of the well-known clover-leaf structure of tRNA is
$\ob \ob \cb \ob \cb \ob \cb \cb$. The asymptotic number of shapes for a length $n$
sequence yields the run time for the Giegerich Lab software {\tt RNAshapes} on length $n$
sequences, since Steffen et al. \cite{giegerich:shapesBioinf} report that
{\tt RNAshapes} runs in time $O(n^3 k s)$ for $s$ sequences, each of 
length at most $n$ and $k$ shapes.}
In  \cite{hofacker99a} Hofacker et al. showed that the 
asymptotic number of {\em canonical} secondary structures
(those having no isolated base pair) is
$2.1614 \cdot n^{-3/2} \cdot 1.96798^n$, a result that was confirmed
by a different method in Clote et al.  \cite{Clote.jbcb09}, where additionally
the expected number of base pairs was shown to be $0.31724 \cdot n$.

A {\em locally optimal}, or {\em kinetically trapped}, 
secondary structure $S$ is one for which no secondary structure $T$,
obtained from $S$ by the removal or addition of a single base pair, has lower energy. It follows
that saturated structures are exactly the kinetically trapped structures with respect to the
{\em Nussinov energy model} \cite{nussinovJacobson}, in which each base pair receives a 
stabilizing energy contribution of $-1$. In this paper, we consider the 
{\em base stacking energy} model, in which each stacked base
pair receives a stabilizing energy contribution of $-1$. 
Here, the base pair $(i,j)$
in secondary structure $S$ is defined to be a stacked base pair, provided that
$(i-1,j+1)$ is also a base pair in $S$ -- 
i.e. provided that there is an outer base pair that provides a 
stabilizing stacking energy. 
In \cite{eversGiegerich}, Evers and Giegerich describe a dynamic programming
algorithm to enumerate all structures that are locally optimal
with respect to the base stacking
energy model; i.e. those structures in which no stem can be extended. 
The authors
called such structures ``saturated''. 
When a strictly positive minimal value is specified for the
length of every stem, a structure is 
saturated in the sense of Zuker \cite{zuker:1986} if and only if it is
saturated in the sense of Evers and Giegerich \cite{eversGiegerich}. 
However, as mentioned in \cite{Clote.jcb06}, when the lengths of stems 
are not constrained, there are structures that are saturated in the sense of
Evers and Giegerich \cite{eversGiegerich}, but which are not saturated in the sense of
Zuker \cite{zuker:1986}.  For clarity of
exposition, we will call a secondary structure G-saturated if no stem can be extended. In this
paper, we give an enumerative framework based on weighted plane trees that allows us to
enumerate G-saturated structures (as well as recover the enumeration of secondary structures
and of saturated structures). We also consider analogous problems for structures with annotated single-stranded, stacked nucleotides (also called dangles).

\subsection*{Outline of paper}
The plan of this paper is as follows. In Section~\ref{section:definitions},
we define the notions of secondary structure and context free grammar, and
provide context free grammars for various classes of secondary structures
considered in the paper. In that section, we show that the asymptotic number
of secondary structures with annotated dangles, as computed in the partition
function of the Markham-Zuker software {\tt UNAFOLD} \cite{Markham.mmb08},
is $0.63998 \cdot n^{-3/2} \cdot  3.06039^n$, exponentially larger than the
number of all secondary structures
$1.104366 \cdot n^{-3/2} \cdot 2.618034^n$, previously established by
Stein and Waterman \cite{steinWaterman}. This new result provides a partial
explanation for M. 
Zuker's observation (personal communication) that {\tt UNAFOLD}
requires substantially more computation time when dangles are included.\footnote{To the best of our knowledge, {\tt UNAFOLD} is currently the only software that computes the
partition function over all secondary structures in a mathematically rigorous manner.}
In Section~\ref{section:computationalResults}, we describe the computation
of secondary structure melting curves with respect to the Nussinov energy
model and the base stacking energy model. Figure~\ref{fig:meltingCurve} shows that folding is more cooperative in the base stacking energy
model. 
In Section~\ref{section:eric1}, we describe the correspondence between
RNA secondary structures and plane trees, and then give generating
functions for the number of secondary structures and locally optimal
secondary structures, with respect to the Nussinov model and the 
base stacking energy model. In Section~\ref{section:eric2}, we give
asymptotic results on the number of secondary structures and locally
optimal secondary structures, as well as their expected number of base
pairs.
In Section~\ref{section:eric3}, we give similar asymptotic results when
annotations for external dangles are included for each type of structure.
Finally Section~\ref{section:discussion} summarizes our main contributions.

\section{Definitions}
\label{section:definitions}

\begin{definition}[Secondary structure]
\label{def:secStr}
An RNA secondary structure for a given RNA sequence $a_1,\ldots,a_n$
of length $n$ is defined to be a set $S$ of ordered pairs
$(i,j)$, with $1 \leq i<j \leq n$, such that
the following conditions are satisfied.
\begin{description}
\item[]
1. {\em Watson-Crick and wobble pairs:}
If $(i,j) \in S$, then 
$\{ a_i,a_j \} \in \{ \{ A,U \} \{ G,C \} \{ G,U \} \}$.
\item[]
2. {\em No base triples:}
If $(i,j)$ and $(i,k)$ belong to $S$, then $j=k$;
if $(i,j)$ and $(k,j)$ belong to $S$, then $i=k$.
\item[]
3. {\em Nonexistence of pseudoknots:}
If $(i,j)$ and $(k,\ell)$ belong to $S$, then it is not the case that
$i<k<j<\ell$.
\item[]
4. {\em Threshold requirement for hairpins:}
If $(i,j)$ belongs to $S$, then $j-i > \theta$, for a fixed value $\theta\geq 0$; i.e. there must be
at least $\theta$ unpaired bases in a hairpin loop.
\end{description}
\end{definition}
For software, such as {\mfold} \cite{zuker:mfoldWebserver} and {\rnafold}
\cite{hofacker:ViennaWebServer}, to predict
RNA secondary structure, $\theta$ is taken to be $3$; i.e., for reasons
related to steric constraints, every hairpin is required to contain at
least three unpaired bases.

A base pair $(i,j)\in S$ is called a \emph{link}. An element $i$ is said to be \emph{linked} if it is involved
in a link and \emph{free} otherwise. 
A link $(i,j)$ is said
to be \emph{stacked} onto another link $(i',j')$ if $i'=i+1$ and $j'=j-1$.
A \emph{stem} is a maximal sequence $\ell_0,\ldots,\ell_{k}$ of links such that $\ell_{i}$
is stacked onto $\ell_{i+1}$ for $0\leq i\leq k-1$; the value $k$ is called the \emph{length} of the stem.  
In some applications a threshold condition on stems is required:

\begin{description}
\item[] 5. \emph{Threshold requirement for stems:} 
Each stem has length at least $\tau$, 
for a fixed value $\tau\geq 0$. 
\end{description} 
Note that Condition~(5) is of no effect for $\tau=0$. 

In this paper, we are concerned with the asymptotic number of 
locally optimal structures. In order to employ generating functions,  we
will need to assume the homopolymer model (following a convention established
by Stein and Waterman \cite{steinWaterman}), meaning that any position
can pair with any other position (arbitrary base pairs, not only
Watson-Crick and wobble pairs). We thus define a secondary
structure of a {\em homopolymer} of length $n$ to be a set $S$ of
base pairs $(i,j)$, where $1 \leq i < j \leq n$, such that the
previous conditions (2,3,4,5) are satisfied. 

The following notion of context free grammar is used for two reasons:
{\em (1)} to provide a clean and
succinct definition for RNA secondary structure, with respect to a particular
energy model, and {\em (2)} for certain enumeration results.
See Lorenz et al.  \cite{Lorenz.jcb08} for more on context free grammars 
and their application to combinatorics. In particular, we refer the reader
to \cite{Lorenz.jcb08} for an explanation of the DSV method used in this
article, which allows us to go directly from a context free grammar to
a functional equation for generating functions.

\begin{definition}[Context free grammar]
\label{ChomskyNormalForm}
A context free grammar is given by
$G = (V,\Sigma,R,S)$, where $V$ is a finite set of
nonterminal symbols (also called variables), $\Sigma$ is
a disjoint finite set of terminal symbols,
$S \in V$ is the {\em start} nonterminal, and 
\[
R \subset V \times (V \cup \Sigma)^*
\]
is a finite set of production rules. Elements of
$R$ are usually denoted by $A \rightarrow w$, rather than
$(A,w)$. 

If $x,y \in (V \cup \Sigma)^*$ and $A \rightarrow w$ is
a rule, then by replacing the occurrence of $A$ in
$x A y$ we obtain $x w y$. Such a derivation in one
step is denoted by
$xAy \Rightarrow_G xwy$, while the
reflexive, transitive closure
of $\Rightarrow_G$ is denoted $\Rightarrow^*_G$.
The language generated by context free grammar $G$ is
denoted by $L(G)$, and defined by
$$L(G) = \{ w \in \Sigma^* : S \Rightarrow^*_G w \}.$$
\end{definition}

Now, in the following sections, we give context free grammars for
RNA secondary structures, including structures with explicitly annotated
dangles. Using the correspondence between grammar and recursions for dynamic
programming, each grammar corresponds to an algorithm for the partition 
function for secondary structures with respect to a different energy model --
the Nussinov model, the base stacking energy model, the Turner model, the
Turner model with a rigorous treatment of dangles, the Turner model with
external dangles.
For notational simplicity, we take $\theta$, the minimum number of
unpaired bases in a hairpin loop to be $1$ (see condition 4 of
Definition \ref{def:secStr}). It is not difficult to extend
the grammar for any fixed value of $\theta$.\footnote{This is done, for 
instance, in grammar $G_4$ by replacing the 
rule $\ub_{\geq \theta} \rightarrow \ub$ by
$\ub_{\geq \theta} \rightarrow \ub^{\theta}$, where
$\ub^{\theta}$ consists of $\theta$ occurrences of $\ub$.}

\subsection*{Nussinov energy model}

In \cite{nussinovJacobson}, Nussinov and Jacobson describe a dynamic
programming algorithm to compute the minimum energy structure for 
a simple energy model, in which each base pair constitutes an energy term
of $-1$.

It is well-known \cite{Lorenz.jcb08} 
that  the following unambiguous grammar  $G_1$
generates all secondary structures of the homopolymer model with $\theta=1$.
Here $G_1$ has start non-terminal symbol $S$, and terminal symbols
$\bullet, \op, \cp$. The
non-terminal symbol $S$ generates all non-empty secondary structures
by using the following grammar (or production) rules.\footnote{Our grammar $G_1$ is
equivalent to the ``tree grammar {\tt nussinov78}'' from 
\cite{Steffen.bb05}.}
\begin{eqnarray*}
S &\rightarrow & \bullet | \op S \cp | S \op S \cp 
\end{eqnarray*}
Let $S(z)$ denote the complex generating function
$S(z) = \sum_{n=0}^{\infty} s_{n} z^n$, where 
Taylor coefficient
$[z^n]S(z)$ is the number $s_n$ of secondary
structures for a homopolymer of size $n$.
By the DSV methodology \cite{Clote.jbcb09,Lorenz.jcb08}, we have
\begin{eqnarray*}
S(z) = S = z + zS + z^2S + z^2 S^2.
\end{eqnarray*}
Introducing the auxilliary variable $u$ to count number of base pairs, we have
\begin{eqnarray}
\label{eqn:waterman}
S(z,u) = S = z + zS + uz^2S + uz^2 S^2 = \displaystyle
\sum_{n} \sum_{k \leq n} s_{k,n} u^k z^n
\end{eqnarray}
where $s_{k,n}$ denotes the number of secondary structures on a length 
$n$ homopolymer, having exactly $k$ base pairs. It follows that
\begin{eqnarray*}
\frac{\partial S(z,u)}{\partial u} &=& \sum_n \sum_{k\leq n}
 k s_{k,n} u^{k-1} z^n\\
\end{eqnarray*}
hence 
\begin{eqnarray}
\label{eqn:A1}
[z^n] \frac{\partial S(z,u)}{\partial u} (z,1) =
\sum_{k \leq n} k s_{k,n}. 
\end{eqnarray}
Since
\begin{eqnarray*}
[z^n] S(z,1) = \sum_{k \leq n} s_{k,n} \\
\end{eqnarray*}
is the number of secondary structures on a homopolymer of length $n$, it
follows that
the asymptotic expected energy over all secondary structures of a homopolymer
of length $n$,  with respect to the Nussinov energy model,  is
equal to $-1$ times the asymptotic expected number of base pairs
\begin{eqnarray}
\label{eqn:A2}
- \lim_{n \rightarrow \infty} \frac{[z^n] \frac{ \partial S(z,u)}{\partial u} (z,1)}{
[z^n] S(z,1)}. 
\end{eqnarray}

\subsection*{Base stacking energy model}
In the base stacking energy model, an energy 
term of $-1$ is assigned to each base pair $(i,j)$ of structure $S$,
provided that $(i,j)$ has an outer stacking pair -- i.e. provided that
$(i+1,j-1) \in S$.  The set of all secondary structures is
generated by the context free grammar
$G_2$ with non-terminals $S,T$, start symbol $S$, and terminals
$\bullet, \op, \cp$ with the following rules:
\begin{eqnarray*}
S &\rightarrow & \bullet | S \bullet | T | ST \\
T &\rightarrow & \op \bullet \cp | \op S \bullet \cp  | \op T \cp  
| \op S T \cp 
\end{eqnarray*}
Here, the non-terminal $S$ generates all secondary structures, while the
non-terminal $T$ generates all secondary structures, such that the
first and last positions are base-paired together. 
By introducing auxilliary non-terminal $T$, 
we can count the number 
of {\em stacked base pairs}, as well as the number of {\em base pairs}.
It is not difficult to show by induction that $G_2$ is an unambiguous grammar
that generates all secondary structures, hence is equivalent to the
previous grammar $G_1$.

By the DSV methodology \cite{Clote.jbcb09,Lorenz.jcb08}, the generating function
$S(z) = \sum_{n} s_n z^n$ satifies the following equations
\begin{eqnarray*}
S(z) = S  &= &z +zS +T + ST\\
T(z) = T  &= &z^2 T +z^2 ST  +z^3 + z^3 S.
\end{eqnarray*}
Introducing the auxilliary variables $u,v$ responsible for counting
the number of base pairs resp. number of stacked base pairs, we have
\begin{eqnarray}
\label{eqn:watermanBis}
S(z,u,v) = S  &= &z +zS +T + ST\\
T(z,u,v) = T  &= &uvz^2 T +uz^2 ST  +uz^3 + uz^3 S. \nonumber
\end{eqnarray}
Letting $s_{k,m,n}$ denote the number of secondary structures on a length 
$n$ homopolymer, having $k$ stacked base pairs and $m$ base pairs, we have
\begin{eqnarray*}
S(z,u,v) &=& \sum_{n} \sum_{k,m \leq n} s_{k,m,n} u^k v^m z^n\\
\frac{\partial}{\partial u} S(z,u,v) &=& \sum_{n} \sum_{k,m \leq n} k 
s_{k,m,n} u^{k-1} v^m z^n\\
\end{eqnarray*}
hence
\begin{eqnarray}
\label{eqn:1}
[z^n] 
\frac{\partial S(z,u,v)}{\partial u} (z,1,1)
 &=& \sum_{k,m \leq n} k 
s_{k,m,n}. 
\end{eqnarray}
Since $S(z,1,1)$ is
is the number of secondary structures on a homopolymer of length $n$, it follows
that the asymptotic expected energy over all secondary structures of a homopolymer
of length $n$,  with respect to the base stacking energy model,  is
equal to $-1$ times the asymptotic expected number of stacked base pairs,
\begin{eqnarray}
\label{eqn:2}
- \lim_{n \rightarrow \infty} 
\frac{ \frac{\partial S(z,u,v)}{\partial u} (z,1,1)}
{[z^n] S(z,1,1)}.
\end{eqnarray}

\subsection*{Grammar for McCaskill algorithm}

All thermodynamics-based RNA secondary structure prediction
algorithms use the Turner nearest neighbor energy model
\cite{turner,xia:RNA}, which contains free energy parameters for
base stacking, single nucleotide dangles, hairpins, bulges,
internal loops and multiloops. These parameters are obtained by a
least squares fit of UV absorption data in optical melting experiments.
For instance, at $37^{\circ}$ C the RNA-RNA stacking free energy of
$\begin{array}{ll}
\mbox{$5'$-{\tt AC}-$3'$}\\
\mbox{$3'$-{\tt UG}-$5'$}\\
\end{array}$
is $-2.24$ kcal/mol and that of
$\begin{array}{ll}
\mbox{$5'$-{\tt CC}-$3'$}\\
\mbox{$3'$-{\tt GG}-$5'$}\\
\end{array}$
is $-3.36$ kcal/mol \cite{xia:RNA}. 
Software such as {\tt mfold} of Zuker \cite{Zuker:1981lr} and
{\tt RNAfold} from Vienna RNA Package \cite{hofacker:FastFolding}
use the Turner energy model, while alternative approaches,
such as {\tt Pfold} \cite{Knudsen:2003lr}
use stochastic context free grammars.

\begin{figure}
\begin{center}
\includegraphics[width=8cm]{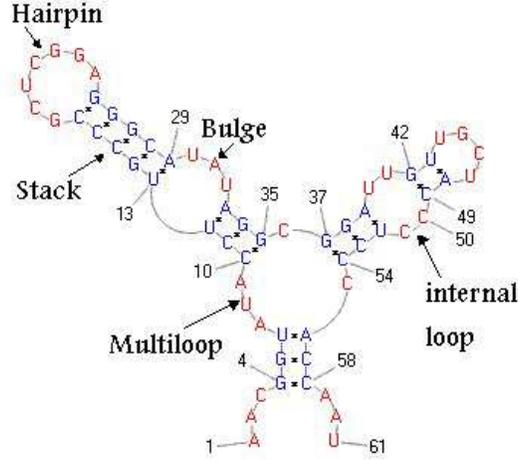}
\end{center}
\caption{Secondary structure together with various loops: stacked base pair,
hairpin, bulge, internal loop, multiloop.}
\label{fig:rnaSecStr2}
\end{figure}

In \cite{mcCaskill}, McCaskill describes a cubic time, dynamic programming
algorithm to compute the partition function $Z = \sum_S \exp(-E(S)/RT)$ over
all secondary structures $S$ of a given RNA sequence. Here $R$ is the universal gas constant,
$T$ is absolute temperature, and $E(S)$ is the energy of structure $S$ with respect to
the Turner energy model \cite{xia:RNA}. By analyzing McCaskill's recursions, we obtain
the following grammar $G_3$, which generates the same set of secondary structures as
that generated by $G_1,G_2$; however, by permitting the classification of various types of
loops, the grammar $G_3$ will permit us later to incorporate energy terms for
{\em dangles}, also known as 
single-stranded, stacked nucleotides, into our considerations.
A {\em stacked base pair} in secondary structure $S$ is given by base pair $(i,j)\in S$, 
such that $(i-1,j+1)\in S$.
A {\em hairpin loop} in secondary structure $S$ is given by base pair $(i,j)\in S$, such that
$i+1,\ldots,j-1$ are unpaired in $S$. A {\em left bulge}
of $S$ is given by base pairs $(i,j),(k,\ell) \in S$, such that $i+1<k<\ell<j$ and 
$\ell+1=j$.  A {\em right bulge}
of $S$ is given by base pairs $(i,j),(k,\ell) \in S$, such that $i<k<\ell<j-1$ and 
$i+1=k$. An {\em internal loop} of $S$ is given by base pairs $(i,j),(k,\ell) \in S$, such that 
$i+1<k<\ell<j-1$; i.e. an internal loop is comprised of both a left and right bulge.
A {\em multiloop} $M$ 
of $S$ is given by base pairs $(i,j),(k_1,\ell_1),\ldots,(k_r,\ell_r) \in S$, 
such that $i<k_1<\ell_1< \cdots < k_r < \ell_r < j$, where $r \geq 2$, and positions
$i+1,\ldots,k_1-1,\ell_1+1,\ldots,k_2-1,\ell_2+1,\ldots,k_r-1,\ell_r+1,\ldots,j-1$ are all
unpaired in $S$. For any positions
$i<x<y<j$, where we do not require $x$ or $y$ to be base-paired, we say that the multiloop
restricted to $[x,y]$ has $c$ components, if exactly $c$ of the base pairs 
$(k_1,\ell_1),\ldots,(k_r,\ell_r)$ are found in the interval $[x,y]$.  See 
Figure~\ref{fig:rnaSecStr2} for an illustration of various loops, and see
\cite{zukerStiegler,zukerMathewsTurner:Guide} for more on loop classification and the Turner
energy model.

Let grammar $G_3$ contain
non-terminal symbols $S$ (start), $U$ (unpaired portion), $B$ (base-paired
portion),
$M_1$ (multiloop with exactly one component),
$M$ (multiloop with at least one component),
with the following production rules
\begin{eqnarray*}
S &\rightarrow & \bullet  | B | S \bullet | SB \\
B &\rightarrow & \op U \cp | \op B \cp |
\op U B \cp | \op B U \cp | \op U B U \cp | \op M M_1 \cp\\
U &\rightarrow & \bullet | \bullet U \\
M_1 &\rightarrow & B | B U\\
M &\rightarrow & M_1 | U M_1 | M M_1\\
\end{eqnarray*}
It is not difficult to show that $G_3$ is an unambiguous context
free grammar equivalent to $G_1,G_2$, thus generates all
secondary structures. The grammar $G_3$ is equivalent to the 
``tree grammar {\tt wuchty98}'' as defined in \cite{Steffen.bb05}, though
notation is vastly different.

\subsection*{Grammar for Markham-Zuker algorithm}

To the best of our knowledge, 
the Markham-Zuker software {\tt UNAFOLD} \cite{Markham.mmb08} is the only 
current 
thermodynamics-based algorithm that computes the partition function for RNA secondary structures
in a mathematically rigorous manner, including correct treatment of 
energy contributions from single-stranded, stacked nucleotides 
-- also called {\em dangles}.  By enlarging the set of terminal symbols, we describe here
an unambiguous context free grammar $G_4$, which generates all secondary structures with
dangle explicitly given. M. Zuker (personal communication) has mentioned that the
algorithm {\tt UNAFOLD} may take approximately twice as long to run when the user chooses to
include treatment of dangles. As we will later see, an explanation for this phenomenon is
that the asymptotic number of secondary structures, where the dangle state is explicitly
annotated, is much larger than the number of secondary structures.

The context free grammar $G_4$ has start symbol $S$, terminal alphabet 
$\{ 5,3,\ub,\op,\cp\}$ and non-terminal alphabet
$\{ S,B,M,M_1,U,\ub_{\geq \theta}  \}$ and rule set
\begin{eqnarray*}
S &\rightarrow &\ub | S \ub | \{ \epsilon + S \} B |
\{ \epsilon+S \} 5 B | \{ \epsilon+S \} B  3 | \{ \epsilon+S \} 5  B  3  \\
B &\rightarrow &\op \ub_{\geq \theta} \cp |
\op B \cp | \op U B \cp | \op B U \cp | \op U B U \cp | \\
& &\op M M_1 \cp | \op 3 M M_1 \cp | \op M M_1 5 \cp | \op 3 M M_1 5 \cp  \\
M &\rightarrow &  \{ \epsilon + U + M\} M_1 \\
M_1 &\rightarrow & M_1 \ub | B | 5 B | B 3 | 5 B 3\\
\ub_{\geq \theta} &\rightarrow & \ub | \ub_{\geq \theta}  U \\
U &\rightarrow & \ub | U \ub
\end{eqnarray*}
Note that $+,\epsilon$ are meta-symbols, used to express the rules
more succinctly. For instance,
$S \rightarrow  \{ \epsilon + S \} B$ is an abbreviation of the rules
$S \rightarrow  B$ and $S \rightarrow  S B$.
The previous rules provide for an unambiguous 
context free grammar that generates all non-empty secondary structures, where
all dangles are explicitly annotated. For instance,
$5 \op \ub \ub \ub \cp$ indicates that in the secondary structure
$\ub \op \ub \ub \ub \cp$, the first position is single-stranded nucleotide which is
$5'$ to the position $2$, and stacks on the base pair $(2,6)$. Similarly,
$\op \ub \ub \ub \cp 3$ indicates that in the secondary structure
$\op \ub \ub \ub \cp \ub$, the last position is single-stranded nucleotide which is
$3'$ to the position $5$ and stacks on the base pair $(1,5)$. 
Since the Turner energy
parameters for hairpins, bulges and internal loops already include contributions for
single-stranded positions within the loop which may dangle on the outer, closing base
pair, it follows that in thermodynamics-based structure prediction, we 
do {\em not} consider  {\em internal} dangles in hairpins,
bulges or internal loops of the 
form $\op 3  \cdots  \cp$, $\op \cdots 5 \cp$, $\op 3 \cdots  5 \cp$, though
such internal dangles are considered in multiloops. Of course, {\em external} dangles
of the form $5 \op \cdots \cp$, $\op \cdots \cp 3$ and
$5 \op \cdots \cp 3$ are considered for all types of loops.

In the grammar $G_4$, non-terminals represent the
following: $S$ denotes the start symbol to generate all structures, 
$B$ indicates that the leftmost and rightmost positions are paired together,
$M$ denotes a substructure located within a  multiloop,
having at least one component (the base pair closing the
multiloop has been generated before non-terminal $M$),
$M_1$ denotes a substructure located within a multiloop,
having exactly one component, where additionally the
leftmost position is paired with a position in the substructure to the right
(though not necessarily the rightmost position).

Note that the Markham-Zuker approach allows dangle annotations of the
rightmost unpaired nucleotide in $\op B B \bullet \cp$ of the
form $\op B B 5 \cp$ or $\op B B 3 \cp$; i.e. where a single-stranded
position occurring between two closing parentheses can be annotated as
either a $5'$-dangle, $3'$-dangle, or no dangle. Indeed,
\begin{eqnarray*}
S \Rightarrow B \Rightarrow \op  M M_1 5 \cp
\Rightarrow \op M_1 M_1 5 \cp
\Rightarrow \op B M_1 5 \cp
\Rightarrow \op B B 5 \cp
\end{eqnarray*}
and
\begin{eqnarray*}
S \Rightarrow B \Rightarrow \op  M M_1  \cp
\Rightarrow \op M_1 M_1 3 \cp
\Rightarrow \op B M_1 3 \cp
\Rightarrow \op B B 3 \cp
\end{eqnarray*}
and
\begin{eqnarray*}
S \Rightarrow B \Rightarrow \op  M M_1  \cp
\Rightarrow \op M_1 M_1 \cp
\Rightarrow \op B M_1 \cp
\Rightarrow \op B M_1 \bullet \cp
\Rightarrow \op B B \bullet \cp
\end{eqnarray*}

\begin{theorem}
In the homopolymer model, where the minimum number of unpaired bases in a hairpin
loop is $1$, the asymptotic number of secondary structures with annotated dangles,
following the Markham-Zuker recursions in \cite{markhamPhD} is
\[
S_n \sim 0.63998 \cdot n^{-3/2} \cdot  3.06039^n .
\]
\end{theorem}
\noindent
{\sc Proof sketch:} It is not difficult to prove by recursion on $n$ that the
set of dangle-annotated secondary structures of length $n$
generated by grammar $G_4$, is equal to the value of the Markham-Zuker
partition function described in pages 14-16 of \cite{markhamPhD}, provided that
all energies are set to $0$.\footnote{It is clear that the number of
structures equals the partition function $\sum_{S} \exp(-E(S)/RT)$ provided that
$E(S)=0$.} Now apply DSV methodology and analyze the dominant singularity
using the Flajolet-Odlyzko theorem, as fully described in \cite{Lorenz.jcb08}.
(At
\url{http://bioinformatics.bc.edu/clotelab/}, 
we provide a detailed computation using Mathematica.)  $\Box$

\subsection*{Grammar for external dangles}

Define {\em external dangle} to mean a $5'$-dangle, which occurs to
the immediate left of an opening parenthesis, or a $3'$-dangle, which
occurs to the right of a closing parenthesis. 
Since our work is theoretical in nature, in the construction
using plane trees in Section~\ref{section:eric3}, we choose to consider the
case that all dangles are external; i.e. no internal dangles,
such as the earlier examples of $\op B B 5 \cp$ and $\op B B 3 \cp$,
are allowed.
We now give a context free grammar for secondary structures having possible
$5'$-dangles and $3'$-dangles in bulges, internal loops, multiloops and external
loops.
Let $G_5$ be a context free grammar with start symbol $S$, 
terminal alphabet $\{ 5,3,\ub,\op,\cp\}$ and non-terminal alphabet
$\{ S,B,M,M_1,U,\ub_{\geq \theta}  \}$ and rule set
\begin{eqnarray*}
S &\rightarrow &\ub | S \ub | \{ \epsilon + S \} B |
\{ \epsilon+S \} 5 B | \{ \epsilon+S \} B  3 | \{ \epsilon+S \} 5  B  3  \\
B &\rightarrow &\op \ub_{\geq \theta} \cp |
\op B \cp | \op \{5 + U5 + U \} B \cp | \op B \{ 3 + 3U + U \} \cp | \\
& &\op \{5 + U5 + U \} B \{ 3 + 3U + U \} \cp | 
\op M M_1 \cp  \\
M &\rightarrow &  \{ \epsilon + U + M\} M_1 \\
M_1 &\rightarrow & M_1 \ub | B | 5 B | B 3 | 5 B 3\\
\ub_{\geq \theta} &\rightarrow & \ub | \ub_{\geq \theta}  U \\
U &\rightarrow & \ub | U \ub
\end{eqnarray*}
As in grammar $G_4$, the symbols $+,\epsilon$ are meta-symbols, to permit
a concise representation of grammar rules; moreover, the meaning of 
non-terminals $S,B,M,M_1$ is the same in $G_5$ as in $G_4$.
It can be proved by induction that grammar $G_5$ is an unambiguous 
context free grammar, that generates all non-empty secondary structures 
with explicitly annotated $5'$-dangles and 
$3'$-dangles, i.e. those dangles that are 
external to any type of loop, whether the loop is
a hairpin, bulge, internal loop, multiloop or external loop.

\begin{theorem}
In the homopolymer model, where the minimum number of unpaired bases in a hairpin
loop is $1$, the asymptotic number of secondary structures with annotated external dangles,
generated by grammar $G_5$ is
\[
S_n \sim 0.96691 \cdot  n^{-3/2} \cdot 3.079596^n. 
\]
\end{theorem}
\noindent
{\sc Proof sketch:} 
Using the DSV methodology, we analyze the dominant singularity
using the Flajolet-Odlyzko theorem, as fully described in \cite{Lorenz.jcb08}.
At  
\url{http://bioinformatics.bc.edu/clotelab/}, 
we provide a detailed computation using Mathematica. 
Moreover, in the latter part of this paper, in a self-contained manner, we give
an alternate proof using plane trees.  $\Box$

\subsection*{Grammar for saturated structures}

In \cite{Clote.jbcb09}, we presented the following grammar which generates all
saturated secondary structures in the sense of Zuker \cite{zuker:1986}; i.e.
locally optimal with respect to the Nussinov energy model.
Let $G_6$ be the context-free grammar with nonterminal symbols $S,R$,
terminal symbols $\bullet, \op,\cp$, start symbol $S$ and production rules
\begin{eqnarray*}
\label{eqn:CFGforSatStrS}
S & \rightarrow & \bullet | \bullet \bullet | R \bullet | R \bullet \bullet | \op S \cp |  S \op S \cp \\
\label{eqn:CFGforSatStrR}
R & \rightarrow & \op S \cp |  R \op S \cp
\end{eqnarray*}
It can be shown by induction on expression length that
$L(S)$ is the set of saturated structures, and $L(R)$ is the
set of saturated structures with no {\em visible} position; i.e.
external to every base pair.
Here, position $i$ is said to be visible in a secondary structure $T$
if it is external to every base pair of $T$; i.e. for all $(x,y) \in T$,
$i<x$ or $i>y$.

It is possible to describe context free grammars that generate
(1) all secondary structures, (2) all saturated secondary structures,
(3) all G-saturated secondary structures, optionally 
with annotated external dangles. However, the subsequent analysis of 
dominant singularity becomes increasingly arduous. For this reason, 
beginning in Section~\ref{section:eric1}, we
present a new, unified method using duality, marked plane trees,
substitution of generating functions, and the Drmota-Lalley-Woods 
theorem (see Theorem~\ref{thm:drmotaLalleyWoods}).

\section{Computational results}
\label{section:computationalResults}

\begin{figure*}
\begin{center}
\includegraphics[width=0.8\textwidth]{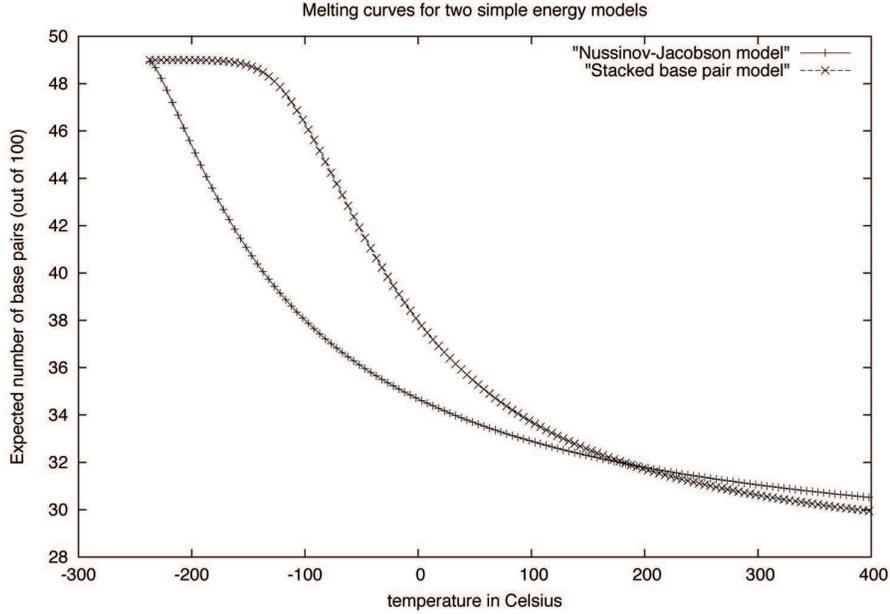}
\caption{Theoretical melting curve for two simple energy models of RNA secondary
structure. Temperature in Celsius is given on the $x$-axis, while expected number
of base pairs is given on the $y$-axis. We implemented an algorithm, 
using dynamic programming, 
with run time $O(n^5)$ and space $O(n^3)$, to compute the
partition function $Z_k = \sum_{S \in \mathbb{S}_k} \exp(-E(S)/RT)$, where
$\mathbb(S)_k$ denotes the set of all secondary structures for a homopolymer
of length 100 nt, having exactly $k$ base pairs. The expected number of
base pairs is thus $\sum_k k \cdot p_k$, where 
$p_k = \frac{Z_k}{Z}$ denotes the probability that a
secondary structure has $k$ base pairs, and $Z$
denotes the full partition function $Z = \sum_{S} \exp(-E(S)/RT) = \sum_k Z_k$.
(Alternatively, and more simply,
we could have produced this curve from the Taylor coefficients of the
expressions to the right of the limit in equations~(\ref{eqn:A2}) and(\ref{eqn:2}), after first solving for $S(z,u)$ 
[resp. $S(z,u,v)$] in equation~(\ref{eqn:waterman}) [resp. 
(\ref{eqn:watermanBis})].)
In the Nussinov-Jacobson energy model \cite{nussinovJacobson}, $E(S)$ is
defined to be $-1\cdot |S|$; i.e. $-1$ times the number of base pairs of $S$.
In the base stacking energy model, $E(S)$ is
defined to be $-1$ times the number of {\em stacked} base pairs of $S$.
Although both models are quite similar, we see that the melting curves are
indeed different, where the base stacking model entails more {\em cooperative}
folding (see \cite{dillBromberg} for discussion of cooperative folding).  }
\label{fig:meltingCurve}
\end{center}
\end{figure*}

In this section, we present computational results to highlight differences
between the Nussinov model and the base stacking energy model, and additionally
to determine the relation between folding time and number of saturated 
structures.
Figure~\ref{fig:meltingCurve}
displays a {\em melting curve} with respect to the Nussinov energy model
and the base stacking energy model. By extending ideas we first described in
\cite{Clote.jcb05}, we developed two algorithms (one for the
Nussinov model and one for the base stacking energy model), each
running in time $O(n^5)$ and space $O(n^3)$, to compute the {\em expected number} of
base pairs as a function of temperature.\footnote{
Alternatively, and more simply, 
we could have produced this curve from the Taylor coefficients of the 
expressions to the right of the limit in equations~(\ref{eqn:A2}) and 
(\ref{eqn:2}), after first solving for $S(z,u)$ 
[resp. $S(z,u,v)$] in equation~(\ref{eqn:waterman}) [resp. 
(\ref{eqn:watermanBis})].} 
Figure~\ref{fig:meltingCurve} clearly shows
that the {\em melting temperature} $T_M$, depends on the energy model, where $T_M$ is
defined as the temperature at which, on average, half the base pairs of the high temperature structure are 
no longer present. Moreover, as the figure shows, the base stacking energy model leads to
more {\em cooperative} folding, as signified by the sigmoidal nature of the curve
(see Dill and Bromberg \cite{dillBromberg} for a discussion of cooperative folding). 

Additionally, the Nussinov energy model and the base stacking energy model are remarkably 
different with respect to pseudoknotted structures, defined by dropping requirement (3) in our
definition of secondary structure; i.e. a pseudoknotted structure $S$ allows base pair 
crossings of the form $(i,j), (k,\ell) \in S$, where $i<k<j<\ell$. While Tabaska et al.
\cite{tabaskaCaryGabowStormo} showed that the minimum energy pseudoknotted structure
can be computed in cubic time $O(n^3)$ by using the maximum weighted matching algorithm,
provided one considers the Nussinov energy model, in the
preprint \cite{ponty:pseudoknotNPcomplete}, Sheikh et al.
show that determination of the minimum energy pseudoknotted structure
for the base stacking energy model is $NP$-complete, a refinement of a result
of Lyngs{\o} and Pedersen \cite{Lyngso.jcb00}.

\section{Enumeration of locally optimal  secondary structures}
\label{section:eric1}

\subsection{Duality: RNA secondary structure $\leftrightarrow$ weighted plane tree}

\begin{figure}[t!]
\begin{center}
\includegraphics[width=15cm]{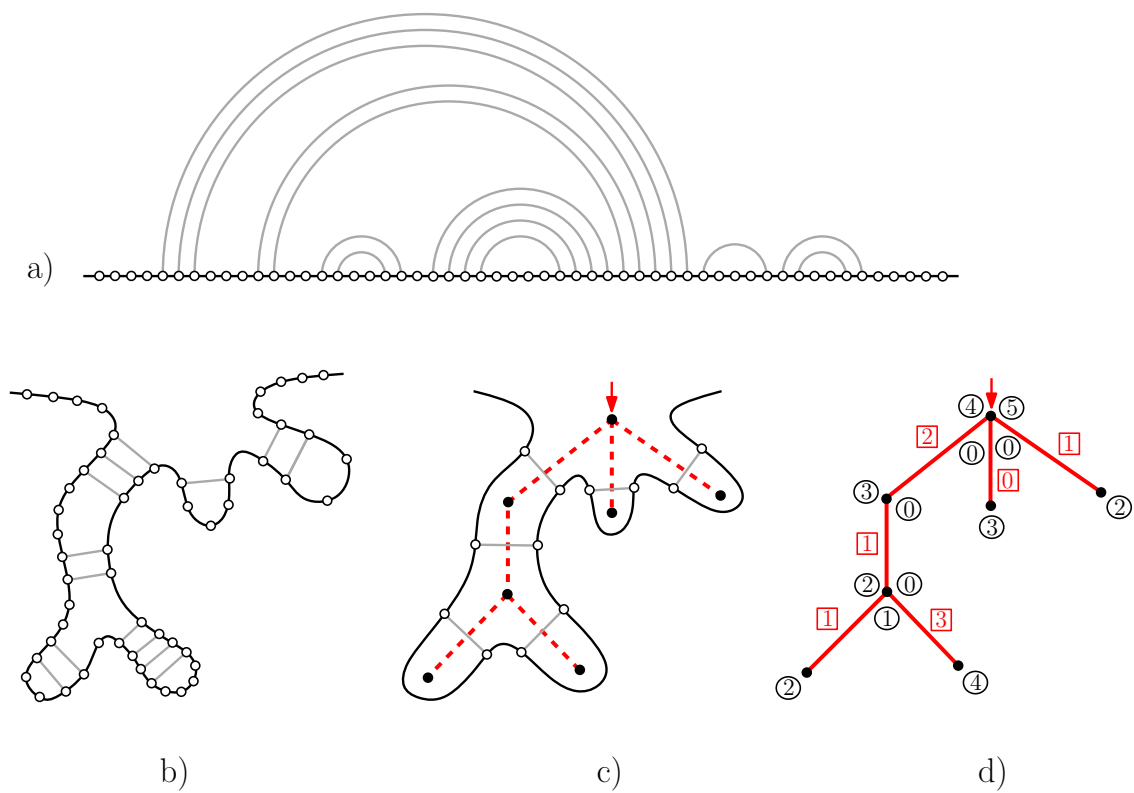}
\end{center}
\caption{(a) A (homopolymer) secondary structure, (b) deformed into a tree-like shape, (c) the reduced structure superimposed
with the dual rooted plane tree (in dashed lines, with the root indicated by an ingoing arrow), (d) the rooted plane tree with weights at corners (surrounded
by circles) 
to indicate segment lengths, and weights at edges (surrounded by squares) to indicate stem lengths.}
\label{fig:RNA_tree}
\end{figure}

It is well known that secondary structures have a tree shape, and there are 
several ways to formulate it. 
Here we find convenient to associate in a bijective way to a
secondary structure (in the homopolymer formulation)  
a rooted plane tree with nonnegative integers (weights) at the corners and at the edges.
The transformation is shown in Figure~\ref{fig:RNA_tree}. 
Start with a secondary structure $S$ of length $n$, the elements in the sequence being ranked from $1$ to $n$.  
Call \emph{segment} of $S$ a sequence $i,i+1,\ldots,j$ such
that $i<j$ and: (i)~either $i=0$, or $1\leq i\leq n$ and the element $i$ is linked,
(ii)~either $j=n+1$, or $1\leq j\leq n$ and the element $j$ is linked, 
(iii)~all elements in $i+1,\ldots,j-1$ are free. 
Note that there are 
$j-i-1$ free elements in the segment.   
Then perform two reduction operations on~$S$:
\begin{description}
\item[Stem-reduction]
Replace each stem $\ell_0,\ldots,\ell_k$ by a single link.
\item[Segment-reduction]
Replace each segment by a unit segment (with no free element on it).
\end{description}

Call $R$ the reduced structure (which has no free element).  
Given the standard plane representation of $R$, 
draw a vertex, called a \emph{dual vertex} in each region, and for each 
link of $R$, draw a \emph{dual edge}
connecting the vertices in the regions on each side of the link.  
 The obtained figure (keeping the dual vertices and dual edges only) is a rooted plane tree $T$.
Note that each edge of $T$ corresponds to a link of $R$ (hence corresponds to a stem of $S$), and each   
corner of $T$ 
 corresponds to a segment of $S$. We \emph{weight} $T$ by giving to each of its corners
a weight corresponding to the number of free elements in the corresponding segment, and giving to each of its
edges a weight corresponding to the length of the corresponding stem.  
Several parameters are in correspondence through the bijection (we use the standard
terminology for parameters of secondary structures, a node of a tree is called a \emph{leaf} if its arity is $0$ 
and an \emph{inner node} if it has positive arity): 
See Table~\ref{table:correspondenceStrTree} for a summary of the 
correspondences between secondary structure loops and nodes of a weighted 
tree.

%

\begin{table}[t!]
\begin{center}
\begin{tabular}{rcl}
secondary structure $S$ & $\leftrightarrow$ & weighted tree $T$ \\
\hline
hairpin & $\leftrightarrow$ & leaf \\
bulge & $\leftrightarrow$ & inner node with one child \\
multiloop & $\leftrightarrow$ & inner node with several children \\
segment with $L$ free elements & $\leftrightarrow$ &   corner of weight $L$\\
stem of length $k$ & $\leftrightarrow$ & edge of weight $k$
\end{tabular}
\end{center}
\caption{Correspondence between types of loop in secondary structure $S$ and 
types of node in the plane tree $T$ obtained by duality.}
\label{table:correspondenceStrTree}
\end{table}

Note also that the number of links of $S$ is the number $|E|$ of edges plus the 
total weight $W_e$ over all edges, 
and that the number of free elements of $S$ is the total weight $W_c$ over all corners, 
hence the length $n$ of $S$ satisfies $n=2|E|+2W_e+W_c$.

\vspace{.4cm}

A weighted rooted plane tree with at least one edge is called \emph{admissible} if it corresponds to a valid  secondary structure (which has at least one link), 
i.e., if the weights satisfy the following conditions:
\begin{description}
\item[]
1. Each non-root node with one child has at least one of its two incident corners of positive weight
 (otherwise the stem-reduction would not have been complete).
\item[]
2. Each corner at a leaf has weight at least $\theta$ (to satisfy the $\theta$-threshold condition).
\item[]
3. Each edge has weight at least $\tau$ (to satisfy the $\tau$-threshold condition). 
\end{description}

\subsection{Generating functions}\label{sec:GF}

For $r\geq 1$, a \emph{weighted combinatorial class indexed by $r$ parameters} is a set $\cA$
together with a \emph{weight-function} $W$ from $\cA$ to $\mathbb{R}$ and
 $r$ \emph{parameter-functions} $P_1,\ldots,P_r$ (one for each parameter) from $\cA$ to $\NN$
such that for any fixed integers $n_1,\ldots,n_r$, the set 
 of structures $\gamma\in\cA$
such that $P_1(\gamma)=n_1,\ldots,P_r(\gamma)=n_r$ is finite. This set is denoted $\cA[n_1,\ldots,n_r]$. 
The corresponding multivariate generating function is 
\begin{equation}
A(x_1,\ldots,x_r):=\sum_{\gamma\in\cA}x_1^{P_1(\gamma)}\cdots x_r^{P_r(\gamma)}W(\gamma).
\end{equation}
We say that variable $x_i$ \emph{marks} the parameter $P_i$, for $1\leq i\leq r$.  
We also use the notation 
$$
[x_1^{n_1}\ldots x_r^{n_r}]A(x_1,\ldots,x_r):=\sum_{\gamma\in\cA[n_1,\ldots,n_r]}W(\gamma).
$$ 
In general we consider \emph{enumerative} generating functions, where $W(\cdot)$ assigns weight $1$
to each structure. However we allow ourselves to weight these structures, e.g., to weight
each secondary structure by $p^{\#(\mathrm{links})}$, 
with $p$ a so-called \emph{stickiness parameter}.
The variables $x_i$ are a priori considered as formal, but one can also evaluate
a generating function at given values, provided the sum converges. The \emph{convergence
domain} of $A(x_1,\ldots,x_r)$ is the set of $r$-tuples
 $(x_1,\ldots,x_r)$ of nonnegative real values such that 
 $A(x_1,\ldots,x_r)$ converges.

As a first example, we briefly recall here how to enumerate (homopolymer) secondary structures, 
via the dual representation by weighted rooted plane trees and using generating functions.
Let $\cF$ be the family of rooted plane trees, possibly reduced to a single vertex, 
with some marked corners (to be occupied by positive weights later on) incident to inner nodes such that each node
with one child has at least one marked corner. Let $F\equiv F(u,v,x)$ be the generating function of $\cF$  
 where $u$ marks the number of leaves, 
$v$ marks the number of marked corners, 
and $x$ marks the number of edges. 
When the root-node $v$ has arity $1$, exactly one of its two corners 
is marked, hence the generating function for trees in $\cF$ whose root-node has arity $1$ is $2vxF$.  
When the root-node $v$ has arity $k\geq 2$, there are $(k+1)$ corners
incident to $v$, and each of these can be marked (independently). Hence the generating
function for trees in $\cF$ where the root-node has arity $k$ is 
$(1+v)^{k+1}x^kF^k$.  Consequently, $F$ satisfies
\begin{equation}
F=u+(2v+v^2)xF+\sum_{k\geq 2}x^k(1+v)^{k+1}F^k=u+\frac{x(1+v)^2F}{1-x(1+v)F}-xF. 
\end{equation}
Let $\cG$ be the family of rooted plane trees with at least one edge and 
with some marked corners (to be occupied by positive weights later on) incident to inner nodes such that each non-root node
with one child has at least one marked corner. Let $G\equiv G(u,v,x)$ be the generating function of $\cG$  
 where $u$ marks the number of leaves, 
$v$ marks the number of marked corners, and $x$ marks the number of edges. 
Again by decomposing at the root, we get
\begin{equation}
G=\sum_{k\geq 1}x^k(1+v)^{k+1}F^k=\frac{x(1+v)^2F}{1-x(1+v)F}.
\end{equation}
Let $g(t,s)$ be the series counting secondary structures with at least one link, 
where $t$ marks the number of free elements, and $s$ marks the number of links. Note that
$g(t,s)$ is also the generating function of admissible rooted weighted plane trees where $t$ marks
the total weight over corners, and $s$ marks the number of edges plus the total weight over edges. 
Such a tree is uniquely obtained from a tree in $\cG$ where each corner at a leaf is assigned a weight of value at least $\theta$,
each non-marked corner at an inner node is assigned weight $0$, each marked corner is assigned a positive weight,
and each edge is assigned a weight of value at least $\tau$. Hence we have $g(t,s)=G(U,V,X)$, where
$$
U:=\sum_{i\geq \theta}t^i=\frac{t^{\theta}}{1-t},\ V=\frac{t}{1-t},\ X:=s\sum_{i\geq \tau}s^i=\frac{s^{\tau+1}}{1-s}. 
$$

To summarize, we have an expression (written as a system of two equations) for  the generating function $g(t,s)$ enumerating
secondary structures with at least one link, where $t$ marks the number of free elements and $s$ marks the number of links
(the generating function of all secondary structures, including the ones with no link, is clearly $g(t,s)+t+t^2+\cdots=g(t,s)+\frac{t}{1-t}$).  
Indeed, if we define $f(t,s):=F(U,V,X)$, then we easily see (since the substitutions
of variables are rational expressions whose series-expansion have nonnegative coefficients) 
that there is a one-line equation specifying $f(t,s)$,
of the form $f(t,s)=Q(t,s,f(t,s))$, 
with $Q\equiv Q(t,s,y)$ a rational expression whose series-expansion (in $s$, $t$, $y$)
has nonnegative coefficients. And there is a rational expression $R\equiv R(t,s,y)$ whose series-expansion
has nonnegative coefficients and such that $g(t,s)=R(t,s,f(t,s))$. Precisely
$$
Q=\mathrm{substitute}\left(u=\frac{t^{\theta}}{1-t},v=\frac{t}{1-t},x=\frac{s^{\tau+1}}{1-s}\right)\ \ \mathrm{into}\ \ u+\frac{x(1+v)^2y}{1-x(1+v)y}-xy,
$$
$$
R=\mathrm{substitute}\left(v=\frac{t}{1-t},x=\frac{s^{\tau+1}}{1-s}\right)\ \ \mathrm{into}\ \ \frac{x(1+v)^2y}{1-x(1+v)y}.
$$ 
This allows us to extract the counting coefficients. Let $g_p(t)$ be the weighted generating
function of secondary structures where $t$ marks the length, and where each structure has
weight $p^{\#(\mathrm{links})}$: $g_p(t)=g(t,pt^2)+t/(1-t)$ (the term $t/(1-t)$
gathers secondary structures with no link); for instance for $\theta=1$ and $\tau=0$ we find
$$
g_p(t)=t+t^2+(1+p)t^3+(1+3p)t^4+(1+6p+p^2)t^5+(1+10p+6p^2)t^6+(1+15p+20p^2+p^3)t^7+\cdots.
$$

\subsection{Counting saturated structures}\label{sec:sat}

The Nussinov energy $E(S)$ of a secondary structure $S$ is defined as $E(S)=-L(S)$,
with $L(S)$ the number of links in $S$. 
A secondary structure $S$ is called \emph{saturated} (or \emph{locally optimal} for the Nussinov energy)
if it is not possible to add a link to $S$ (i.e., decrease the energy by $1$) while keeping a valid
 secondary structure.

\begin{figure}
\begin{center}
\includegraphics[width=10cm]{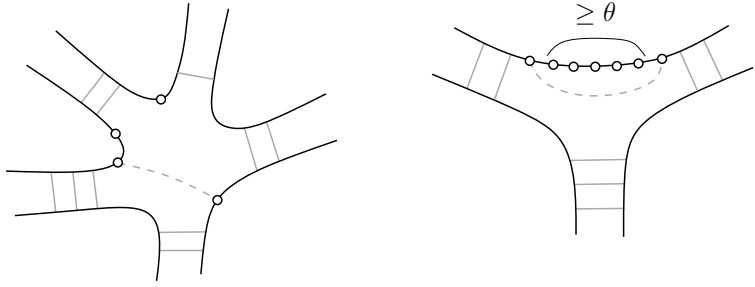}
\end{center}
\caption{Situations where it is possible to add a link to a secondary structure.}
\label{fig:sat}
\end{figure}

\begin{lemma}\label{lem:charact_sat}
Assume $\tau=0$ (no restriction on the lengths of stems). 
Saturated secondary structures with at least one link correspond to admissible weighted rooted plane trees such that:   
\begin{itemize}
\item all corners have weight at most $\theta+1$,
\item at each node there is at most one corner of strictly positive weight.
\end{itemize}
\end{lemma}
\proof
As shown in Figure~\ref{fig:sat}, if there are two positive corners at the same inner node,  
then it is possible to add a link. Also, if there is a corner with weight at least
$\theta+2$ then one can link the first and last free elements in the 
corresponding segment. Hence the weight of each corner is at most $\theta+1$.  
And these are the only two situations where it is possible to add a link
without breaking planarity nor breaking the $\theta$-threshold condition. 
\hfill $\Box$ \medskip

Call $\cF$ the family of rooted plane trees with some marked corners incident to inner nodes (these
marked corners are to be occupied by positive weights later on)  such that:  
(i)~each node with one child has exactly one marked corner, (ii)~each node with several children
has at most one marked corner. Let $F\equiv F(u,v,x)$ be the generating function of $\cF$
where $u$ marks the number of leaves, $v$ marks the number of marked corners, and $x$
marks the number of edges. When the root-node $v$ has arity $1$, exactly one of its two corners 
is marked, hence the generating function for trees in $\cF$ whose root-node has arity $1$ is $2vxF$.  
When the root-vertex $v$ has arity $k\geq 2$, there are $(k+1)$ corners
incident to $v$, and at most one of these corners has positive weight. Hence the generating
function for trees in $\cF$ where the root-vertex has arity $k$ is 
$(1+(k+1)v)x^kF^k$.  Consequently, $F$ satisfies
$$
F=u+2vxF+\sum_{k\geq 2}(1+(k+1)v)x^kF^k,
$$ 
Hence, using the identity $\sum_{k\geq 0}(k+1)A^k=1/(1-A)^2$, $F$ satisfies
\begin{equation}\label{eq:Fsat}
F=u+\frac{x^2F^2}{1-xF}+\frac{v}{(1-xF)^2}-v.
\end{equation}

Now let $\cG$ be the family of rooted plane trees with at least one edge, and with
marked corners incident to inner nodes such that: (i')~each node $v$ with one child has exactly one marked
corner if $v$ is different from the root-node, and has \emph{at most} one marked corner if $v$ is the root-node, 
(ii)~each node with several children
has at most one marked corner. 
Let $G\equiv G(u,v,x)$ be the generating function of $\cG$ where $u$, $v$, $x$ mark respectively
the numbers of leaves, marked corners, and edges. Decomposing again at the root, we get
\begin{equation}
G=\sum_{k\geq 1}(1+(k+1)v)x^kF^k=\frac{xF}{1-xF}+\frac{v}{(1-xF)^2}-v.
\end{equation}

We take here $\tau=0$ (no restriction on the lengths of stems). 
Let $g(t,s)$ be the generating function of saturated secondary structures with at least one link, where
$t$ marks the number of free elements and $s$ marks the number of links.
Then Lemma~\ref{lem:charact_sat} ensures that $g(t,s)=G(U,V,X)$, where
$$
U=t^{\theta}(1+t),\ V=t+\ldots+t^{\theta+1}=\frac{t-t^{\theta+2}}{1-t},\ X=\frac{s}{1-s}.
$$

To summarize (in a similar way as for 
general structures), we have an expression (written as a system of two equations) for  the generating function $g(t,s)$ enumerating \emph{saturated} 
secondary structures with at least one link, where $t$ marks the number of free elements and $s$ marks the number of links
(the generating function of all saturated secondary structures, including the ones with no link, is $g(t,s)+t+\cdots+t^{\theta+1}=g(t,s)+\frac{t-t^{\theta+2}}{1-t}$).
Indeed, if we define $f(t,s):=F(U,V,X)$, then  there is a one-line equation specifying $f(t,s)$,
of the form $f(t,s)=Q(t,s,f(t,s))$, 
with $Q(t,s,y)$ a rational expression whose series-expansion (in $s$, $t$, $y$)
has nonnegative coefficients. And there is a rational expression $R(t,s,y)$ whose series-expansion
has nonnegative coefficients and such that $g(t,s)=R(t,s,f(t,s))$. Precisely
$$
Q=\mathrm{substitute}\left(u=t^{\theta}(1+t),v=\frac{t-t^{\theta+2}}{1-t},x=\frac{s}{1-s}\right)\ \ \mathrm{into}\ \ u+\frac{x^2y^2}{1-xy}+\frac{v}{(1-xy)^2}-v,
$$
$$
R=\mathrm{substitute}\left(v=\frac{t-t^{\theta+2}}{1-t},x=\frac{s}{1-s}\right)\ \ \mathrm{into}\ \ \frac{xy}{1-xy}+\frac{v}{(1-xy)^2}-v.
$$
Again this allows us to extract the counting coefficients. Let $g_p(t)$ be the weighted generating
function of saturated secondary structures where $t$ marks the length, and where each structure has
weight $p^{\#(\mathrm{links})}$: $g_p(t)=g(t,pt^2)+t+\cdots+t^{\theta+1}$; for $\theta=1$ and $\tau=0$ we find
$$
g_p(t)=t+t^2+pt^3+3pt^4+(4p+p^2)t^5+(2p+6p^2)t^6+(17p^2+p^3)t^7+\cdots.
$$ 
 
\subsection{Counting G-saturated structures}\label{sec:Gsat}

The \emph{base stacking energy} $E(S)$ of a secondary structure $S$ is defined as $E(S):=-T(S)$, with $T(S)$ the sum of sizes of all stems of $S$.  
A (homopolymer) secondary structure is called \emph{G-saturated} (locally optimal for the base stacking energy)
if it is not possible to add a link so as to extend a stem (i.e., decrease by $1$ the base stacking energy).   
In general, the addition of a link to a secondary structure either
 creates a new stem of length $0$ or extends an already existing stem. 
Hence, in a G-saturated structure a valid link addition always creates a new stem of length $0$.
In case $\tau>0$, creating a stem of length $0$ is not a valid link addition (since
the stems must have positive length), hence no valid link addition to a G-saturated is possible  for $\tau>0$. 
In other words, the concepts of saturated and of G-saturated structures coincide when $\tau>0$ (whereas
for $\tau=0$ the class of saturated structures is strictly contained in the class of G-saturated structures).  
 In this section we enumerate G-saturated structures according to 
the number of free elements and the number of links,  
for any given values of the threshold parameters $\tau$ and $\theta$. 
 
Again we formulate the conditions on the dual representation.
For this purpose we define adjacency of corners. Two corners $c$ and $c'$ of a rooted plane tree $T$ 
are called \emph{adjacent} if they are incident to the same vertex $v$ of $T$ and there is an edge $e$ 
incident to $v$ such that $c$ and $c'$ are the corners incident to $v$ on each side of $e$. 
Note that the two corners on each side of the root (the root is represented as an ingoing arrow in Figure~\ref{fig:RNA_tree}) 
 are considered as adjacent only when the root-node $v$  
has arity $1$ (in which case they are adjacent through the unique edge incident to $v$).   

\begin{figure}
\begin{center}
\includegraphics[width=10cm]{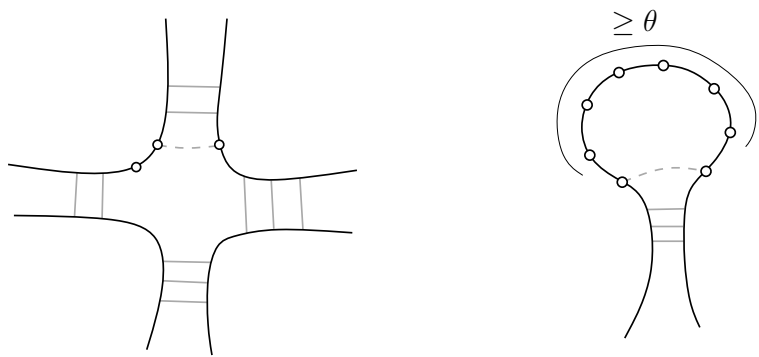}
\end{center}
\caption{Situations where it is possible to extend a stem of a secondary structure.}
\label{fig:Gsat}
\end{figure}
 
\begin{lemma}\label{lem:dualG-sat}
The G-saturated secondary structures with at least one link 
correspond to admissible weighted rooted plane trees such that: 
\begin{itemize}
\item the corners at leaves have weight at most $\theta+1$, 
\item any two adjacent corners can not both have strictly positive weight. 
\end{itemize}
\end{lemma}
\proof
As shown in Figure~\ref{fig:Gsat}, 
 if there are two adjacent positive corners,   
then it is possible to add a link so as to extend an existing stem. 
Also, if there is a corner of weight at least $\theta+2$ at a leaf $\ell$, 
then one can link the first and last free elements in the 
corresponding segment and thus extend the stem associated to the edge leading to $\ell$. 
Hence the weight of a corner at a leaf is at most $\theta+1$.  
And these are the only two situations where it is possible to extend a stem 
without breaking planarity nor breaking the $\theta$-threshold and $\tau$-threshold condition. 
\hfill $\Box$ \medskip

Call $\cF$ the family of rooted plane trees  with
some marked corners incident to inner nodes (again these
marked corners are to be occupied by positive weights later on) such that:  
(i)~each inner node with one child has exactly one marked corner, (ii)~two corners can not both be marked if they 
are adjacent or if they are the two corners on each side of the root (the root is indicated by an ingoing arrow in Figure~\ref{fig:RNA_tree}).     
Let $F\equiv F(u,v,x)$ be the generating function of $\cF$
where $u$ marks the number of leaves, $v$ marks the number of marked corners, and $x$
marks the number of edges. Finding an equation satisfied by $F$ is a little more involved
than for saturated structures. At first we  need a preliminary study on independent
sets (i.e., sets containing only pairwise non-adjacent elements) 
on a $k$-sequence or on a $k$-cycle. 

For $k>0$ and $m\leq k$, let $c_{k,m}$ (resp. $s_{k,m}$) be the number of
ways of choosing $m$ \emph{marked} elements on the oriented cycle  $(1,2,\ldots,k)$
(resp. sequence $1,2,\ldots,k$)  
of $k$ elements such that no two consecutive elements are marked, and let 
$C_k=C_k(v):=\sum_{m}c_{k,m}v^m$ (resp. $S_k=S_k(v):=\sum_{m}s_{k,m}v^m$) 
be the corresponding (polynomial) generating function. The polynomials $S_k$
are well-known to be the \emph{Fibonacci polynomials} and satisfy an easy recurrence which
we briefly recompute here. We take
the convention $S_0=1$. 
Let $k\geq 2$. If an independent set on the $k$-sequence 
 starts with a marked element, then the next element is forbidden
 and the remaining $(k-2)$-sequence might be occupied by any independent set;
 this gives a  contribution $vS_{k-2}$ in $S_k$, where the factor $v$
 takes account of the first element being marked. 
 If an independent set on the $k$-sequence 
 starts with a non-marked element, then the remaining $(k-1)$-sequence might be occupied by any independent set;
 this gives a  contribution $S_{k-1}$ in $S_k$. Therefore 
 $$
 S_k=vS_{k-2}+S_{k-1}\ \ \mathrm{for}\ k\geq 2,\ \ S_0=1,\ S_1=1+v.
 $$
Now define $S\equiv S(v,z):=\sum_{k\geq 0}S_k(v)z^k$. The recurrence on $S_k$
above multiplied by $z^k$ and summed over $k\geq 2$ 
yields $S-S_0-zS_1=vz^2S+z(S-S_0)$. With $S_0=1$ and $S_1=1+v$ we obtain
$$
S=\frac{1+vz}{1-z-vz^2}.
$$
Let us go back to independent sets on the $k$-cycle $(1,\ldots,k)$, for $k\geq 3$. 
In such a set, either $1$ is occupied, in which case the adjacent elements
$2$ and $k$ are unoccupied and the remaining segment $3,\ldots,k-1$
might be occupied by any independent set. This gives contribution $vS_{k-3}$
to $C_k$. If $1$ is unoccupied, then the remaining segment $2,\ldots,k$
might be occupied by any independent set; this gives contribution $S_{k-1}$
to $C_k$. Consequently
$$
C_k=vS_{k-3}+S_{k-1}\ \ \mathrm{for}\ k\geq 3.
$$
If the root-node $v$ of a tree in $\cF$
 has arity $1$ then exactly one of its two incident corners is marked (by definition of $\cF$), thus 
 the generating function of trees
in $\cF$ whose root-node has arity $1$ is $2vxF$; if $v$ has arity $k\geq 2$ then  
the marked corners around $v$ form an independent set
(no two consecutive corners are marked). Thus, for $k\geq 2$, the generating function of trees
in $\cF$ whose root-node has arity $k$ is $C_{k+1}(v)x^kF^k$ (since there are $k+1$
corners incident to the root-node).  
Consequently $F$ satisfies 
\begin{eqnarray*}
F&=&u+2vxF+\sum_{k\geq 2}C_{k+1}(v)x^kF^k\\
&=&u+2vxF+\sum_{k\geq 2}\big[ vS_{k-2}+S_k \big]x^kF^k\\
&=&u+2vxF+vx^2F^2S(v,xF)+\big(S(v,xF)-1-(1+v)xF\big).
\end{eqnarray*}
Using  the rational expression of $S$ above and rearranging, we obtain
\begin{equation}\label{eq:y_G_sat}
F=u+2vxF+\frac{1+2vx^2F^2\cdot(1+vxF)}{1-xF-vx^2F^2}-xF-1.
\end{equation}
Now let $\cG$ be the family of rooted plane trees with at least one edge
and where some corners at inner nodes are marked such that (i) each non-root inner node of arity $1$ 
has exactly one marked corner, (ii)~two adjacent corners can not both be marked. 
And let 
$G\equiv G(u,v,x)$ be the generating function of $\cG$ where $u$ marks the number of leaves,
$v$ marks the number of marked corners, and $x$ marks the number of edges. 
The difference between $\cG$ and $\cF$ is at the root-vertex: in $\cG$ 
the two corners on each side of the root are allowed to be both marked when the root-vertex
has more than one child, and are allowed to be both unmarked when the root-vertex has one child. So we have
$$
G=\sum_{k\geq 1}S_{k+1}(v)x^kF^k.
$$
Using the rational expression of $S$ above and rearranging, we obtain the following expression
for $G$ in terms of $F$:
\begin{equation}
G=\frac{xF(1+2v+(1+v)vxF)}{1-xF-vx^2F^2}.
\end{equation}

Now let $g(t,s)$ be the generating function of G-saturated structures with at least one link, 
where $t$ marks the number of free elements and $s$ marks the number of links. By Lemma~\ref{lem:dualG-sat}, 
\begin{equation}
g(t,s)=G(U,V,X),
\end{equation}
where 
$$
U=t^{\theta}(1+t),\ \ V=\frac{t}{1-t},\ \ X=\frac{s^{\tau+1}}{1-s}. 
$$

The conclusion is similar to the other
two cases (general structures, saturated
structures):  we have an expression (written as a system of two equations) for  the generating function $g(t,s)$ enumerating \emph{G-saturated} 
secondary structures with at least one link, where $t$ marks the number of free elements and $s$ marks the number of links
(the generating function of all G-saturated secondary structures, including the ones with no link, is $g(t,s)+t+t^2+\cdots=g(t,s)+\frac{t}{1-t}$). 
Indeed, if we define $f(t,s):=F(U,V,X)$, then 
 there is a one-line equation specifying $f(t,s)$,
of the form $f(t,s)=Q(t,s,f(t,s))$, 
with $Q(t,s,y)$ a rational expression whose series-expansion (in $s$, $t$, $y$)
has nonnegative coefficients. And there is a rational expression $R(t,s,y)$ whose series-expansion
has nonnegative coefficients and such that $g(t,s)=R(t,s,f(t,s))$. Precisely
$$
Q=\mathrm{substitute}\left(u=t^{\theta}(1+t),v=\frac{t}{1-t},x=\frac{s^{\tau+1}}{1-s}\right)\ \ \mathrm{into}\ \ u+2vxy+\frac{1+2vx^2y^2(1+vxy)}{1-xy-vx^2y^2}-1-xy,
$$
$$
R=\mathrm{substitute}\left(v=\frac{t}{1-t},x=\frac{s^{\tau+1}}{1-s}\right)\ \ \mathrm{into}\ \ \frac{xy(1+2v+(1+v)vxy)}{1-xy-vx^2y^2}.
$$
Again this allows us to extract the counting coefficients. Let $g_p(t)$ be the weighted generating
function of G-saturated secondary structures where $t$ marks the length, and where each structure has
weight $p^{\#(\mathrm{links})}$: $g_p(t)=g(t,pt^2)+t/(1-t)$; for $\theta=1$ and $\tau=0$ we find
$$
g_p(t)=t+t^2+(1+p)t^3+(1+3p)t^4+(1+4p+p^2)t^5+(1+4p+6p^2)t^6+(1+4p+17p^2+p^3)t^7+\cdots.
$$

\section{Asymptotic results}
\label{section:eric2}

\subsection{Asymptotic enumeration}\label{sec:asympt_enum}

We show here that the number of structures of length $n$
follows a universal asymptotic behaviour in $c\ \!\gamma^n n^{-3/2}$ (with $c$ and $\gamma$
explicit positive constants), which is typical of tree-structures.
The proof classically relies on the Drmota-Lalley-Woods theorem~\cite[VII.6]{FlSebook}, which we recall at first.
Consider an equation of the form
\begin{equation}\label{eq:DLL}
a(t)=\Phi(t,a(t)),
\end{equation}
where $\Phi(t,y)$ is a rational expression in $t$ and $y$. Such an equation is called \emph{admissible}
if the following conditions are satisfied:
\begin{itemize}
\item
the rational expression 
$\Phi(t,y)$ has a series-expansion in $t$ and $y$ with nonnegative coefficients, is nonaffine in $y$, 
and satisfies~\footnote{We use the subscript notation for partial derivatives.} $\Phi(0,0)=0$ and $\Phi_y(0,0)=0$,   
\item
the unique generating function $y=a(t)$ solution of~\eqref{eq:DLL} is aperiodic, i.e., can not be written as
$a(t)=t^q\tilde{a}(t^p)$ for some integers $p,q$ with $p\geq 2$. 
\end{itemize}
There is an easy criterion to check the aperiodicity condition: 
it suffices to prove that there is some $n_0$ such that $[t^n]a(t)>0$
for $n\geq n_0$.

\begin{theorem}[Drmota-Lalley-Wood]\label{theo:DLL}
\label{thm:drmotaLalleyWoods}
Let $y=a(t)$ be the generating function that is the unique solution of an admissible equation $y=\Phi(t,y)$.  
Then 
$$
[t^n]a(t)\sim c\ \!\gamma^n\ \!\!n^{-3/2},
$$
where $\gamma=1/t_0$, with $(t_0,y_0)$ the unique pair in the convergence domain of $\Phi(t,y)$ that is  
 solution of the \emph{singularity system}:
$$
y=\Phi(t,y),\ \ \  \Phi_y(t,y)=1;
$$
and where 
$$
c=\sqrt{t_0\Phi_{t}(t_0,y_0)/(2\pi \Phi_{y,y}(t_0,y_0))}. 
$$
Moreover, if $\Psi(t,y)$ is a rational expression not constant in $y$, that has a series-expansion with nonnegative coefficients,  
and such that the convergence domain of $\Psi(t,y)$ is contained in the convergence domain of $\Phi(t,y)$, 
then the coefficients of the generating function $b(t):=\Psi(t,a(t))$ behave as
$$
[t^n]b(t)\sim d\ \!\gamma^n\ \!\!n^{-3/2},
$$
where $d=c\cdot \Psi_y(t_0,y_0)$.
\end{theorem}

\begin{remark}
The Drmota-Lalley-Wood theorem is classically proved (e.g. in~\cite[VII.6]{FlSebook}) for polynomial systems (i.e., for 
$\Phi(t,y)$ a polynomial). But one easily checks that, more generally, if $\Phi(t,y)$ is a bivariate series
that diverges at all its singularities, then the conclusions remain the same.
\end{remark}

From the Drmota-Lalley-Wood theorem we obtain

\begin{proposition}\label{prop:asympt1}
Let $p$ be a fixed  positive real value (stickiness parameter). 
Let $g_p(t)$ be the univariate generating function of general (resp. saturated, resp. G-saturated)
homopolymer secondary structures, where $t$ marks the length of the sequence and where
each  structure has weight $p^{\#(\mathrm{links})}$.    

Then, for any values of the threshold-parameters $\theta$ and $\tau$
($\tau=0$ if one considers saturated structures), 
 there are computable positive constants $c$ and $\gamma$ (depending on $\tau$, $\theta$, $p$, and in which
setting: general, saturated, or G-saturated) such that 
$$
[t^n]g_p(t)\sim c\ \!\gamma^n\ \!\!n^{-3/2}.
$$
\end{proposition}
\proof
Recall that, in each of the three settings (general, saturated, G-saturated), $g(t,s)$ denotes
the generating function of secondary structures with at least one link, 
where $t$ marks the number of free elements and
$s$ marks the number of links. We have seen that, in each of the three settings, 
there are two rational expressions $Q(t,s,y)$ and $R(t,s,y)$ that  have  
nonnegative coefficients (in the series-expansion), and  there is an adjoint generating function $f(t,s)$ such that $f(t,s)=Q(t,s,f(t,s))$
and $g(t,s)=R(t,s,f(t,s))$. In addition, the convergence domain of $Q(t,s,y)$
is clearly the same as the convergence domain of $R(t,s,y)$; for instance, for G-saturated structures,
the convergence domain is the set of nonnegative triples $(t,s,y)$ such that $t<1$, $s<1$, 
and $xy+vx^2y^2<1$, where $v=t/(1-t)$ and $x=s^{\tau+1}/(1-s)$. Note that
in all three settings, $f(0,0)=1$ for $\theta=0$ and $f(0,0)=0$ for $\theta>0$. 
If we set  
$a(t):=f(t,pt^2)-\mathbf{1}_{\theta=0}$  (with $\theta$ the threshold parameter)  
 and $b(t):=g(t,pt^2)$, 
then we are in the conditions of the Drmota-Lalley-Wood theorem, 
with $\Phi(t,y):=Q(t,pt^2,y+\mathbf{1}_{\theta=0})-\mathbf{1}_{\theta=0}$ and 
$\Psi(t,y):=R(t,pt^2,y+\mathbf{1}_{\theta=0})$. The conditions for $\Phi$ and $\Psi$
are readily checked, we show now the 
aperiodicity of $a(t):=f(t,pt^2)$ (proving that the $n$th 
coefficient is strictly positive for $n$ large enough).   
Note that it is enough to consider $p=1$ (the strict positivity of $[t^n]f(t)$ does
not depend on $p>0$).  
In each of the three settings (general, saturated, G-saturated), 
$a(t)$ is the enumerative generating function of some
explicit class of rooted weighted plane trees.
 For instance, for saturated structures, $a(t)$ 
counts admissible rooted weighted plane trees with all corners
of weight at most $\theta+1$, with at most one positive corner per node, 
and where each node of arity $1$ has exactly one positive corner.  
For $i\geq \tau$, consider the weighted rooted plane tree $T_i$ made of one
edge $e$ leading to a leaf $\ell$, with weight $1$ (resp. $0$) at the corner 
to the left (resp. right) of the root, with weight $i$ on $e$ and weight $\theta$ on $\ell$.
And consider the tree $T_i'$ defined exactly as $T_i$ except that $\ell$ has weight $\theta+1$.
Note that $T_i$ contributes to $[t^{2i+\theta+3}]a(t)$ and $T_{i}'$ contributes to $[t^{2i+\theta+4}]a(t)$.
Hence $[t^n]a(t)>0$ for all $n\geq 2\tau +\theta+3$, so $a(t)$ is aperiodic. 
In exactly the same way, $a(t)$ is aperiodic in the general setting and in the G-saturated
setting. 

Theorem~\ref{theo:DLL} ensures that there are $c>0$ and $\gamma>0$ such that $[t^n]g(t,pt^2)\sim c\gamma^n\ \!n^{-3/2}$.
Actually, in the case of general and G-saturated structures, 
we have $\gamma>1$ since (according to Theorem~\ref{theo:DLL}) 
there is some $y_0$ such that 
$(1/\gamma,y_0)$ is in the convergence domain of $\Phi(t,y)$, and since clearly any $(t_0,y_0)$
in the convergence domain of $\Phi(t,y)$ satisfies $t_0<1$ (indeed $Q(t,s,y)$ involves the quantity
$1/(1-t)$, in each of the general and in the G-saturated case). 
The generating function $g_p(t)$ (which includes
also secondary structures with no link, as opposed to $g(t,s)$) satisfies
$g_p(t)=g(t,pt^2)+t/(1-t)$ for secondary and for G-saturated structures, and satisfies
$g_p(t)=g(t,pt^2)+t+\ldots+t^{\theta+1}$ for saturated structures. So the additional term
gathering saturated structures with no link has negligible asymptotic contribution in all cases.
\hfill $\Box$ \medskip

For $p=1$, $g_p(t)$ is the enumerative generating function of homopolymer structures.
Another value of interest is $p=3/8$. 
Indeed, if we want to count RNA secondary structures  (each base is labelled by a letter in $\{A,G,C,U\}$) instead of homopolymers,
this corresponds to giving weight $4$ to each free element (because there are $4$ possible labels) 
and giving weight $6$ to each pair of linked elements (because there are $6$ allowed labellings  
out of $4^2=16$, due to the Watson-Crick and wobble pairs). Therefore the corresponding enumerative generating function is $g(4t,6t^2)$. We have
$$
[t^n]g(4t,6t^2)=4^n[t^n]g(t,3t^2/8)=4^n[t^n]g_{3/8}(t).
$$
In other words, $[t^n]g_{3/8}$ is the \emph{expected number} of RNA secondary structures
with the desired properties (general, saturated, or G-saturated) on a random
sequence of size $n$ (i.e., for a random word in $\{A,G,C,U\}^n$). 
  
Table~\ref{fig:table_asympt} shows the asymptotic behaviour of $[t^n]g_p(t)$ for $p=1$
and $p=3/8$ 
in the three settings. (The methodology to compute $\gamma$ for saturated structures using computer algebra tools
is detailed in~\cite{Clote.jbcb09}.)  

\begin{table}
\begin{center}
\includegraphics[width=14cm]{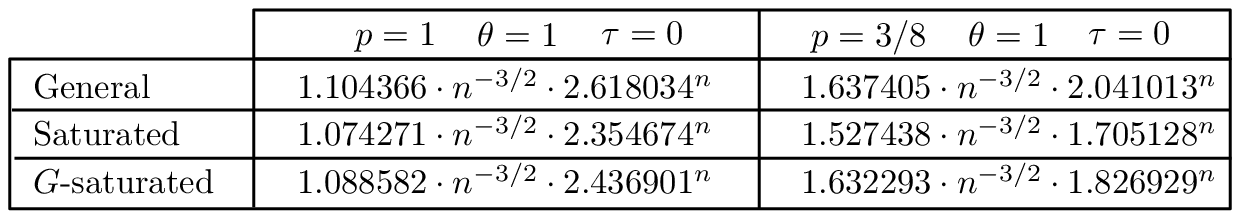}
\end{center}
\caption{Asymptotic behaviour of the $n$th coefficient  of the generating
function $g_p(t)$ counting secondary structures (general, saturated, or G-saturated)
with weight $p$ on each link.}
\label{fig:table_asympt}
\end{table}

\subsection{Limit law for the number of links}
Using a theorem of Drmota~\cite{Dr97} (closely related to the Drmota-Lalley-Wood theorem)
we show that the number of links in a random secondary structure
 (general, saturated, or G-saturated) of length $n$ 
is asymptotically a gaussian law with $\Theta(n)$ expectation and $\Theta(\sqrt{n})$ 
standard deviation.

Consider an equation of the form
\begin{equation}
a(t,u)=\Phi(t,u,a(t,u)),
\end{equation}
where $\Phi(t,u,y)$ is a rational expression in $t$, $u$ and $y$. Such an equation
is called \emph{admissible} if $\Phi(t,u,y)$ is nonconstant in $u$, has a series-expansion (in $t$, $u$, $y$)
with nonnegative coefficients, the equation $y=\Phi(t,1,y)$ is admissible (in the sense
of Section~\ref{sec:asympt_enum}), and there is a $3\times 3$-matrix $m[i,j]$ with integer coefficients and nonzero determinant
such that $[t^{m[i,1]}u^{m[i,2]}y^{m[i,3]}]\Phi(t,u,y)>0$ for all $i\in\{1,2,3\}$. 

\begin{theorem}[Drmota~\cite{Dr97}]\label{theo:drmota}
Let $y=a(t,u)$ be a generating function that is the unique solution of an admissible
equation $y=\Phi(t,u,y)$.  Assume that the generating
function $b(t,u)=\sum_{\gamma\in\cG}t^{|\gamma|}u^{\chi(\gamma)}W(\gamma)$ of a weighted 
combinatorial class $\cG$  
is given by $b(t,u)=\Psi(t,u,a(t,u))$, with $\Psi(t,u,y)$ a rational expression with nonnegative
coefficients (in the series-expansion), nonconstant in $y$,
 and such that the convergence domain of $\Psi(t,1,y)$
 is included in the one of $\Phi(t,1,y)$. For $n\geq 0$ let $\cG_n:=\{\gamma\in\cG,\ |\gamma|=n\}$,
and define the random variable $X_n$ as $\chi(\gamma)$, with $\gamma$ a random structure in $\cG_n$ under the distribution 
$$
P(\gamma)=\frac{W(\gamma)}{\sum_{\gamma\in\cG_n}W(\gamma)}.
$$
For $u>0$ in a neighbourhood of $1$, denote by $\rho(u)$ the radius of convergence of 
$y:t\to a(t,u)$
, and let 
$$
\mu=-\frac{\rho'(1)}{\rho(1)},\ \ \ \sigma^2=-\frac{\rho''(1)}{\rho(1)}-\frac{\rho'(1)}{\rho(1)}+\left(\frac{\rho'(1)}{\rho(1)}\right)^2.
$$
Then $\mu$ and $\sigma$ are strictly positive and
 $\displaystyle\frac{X_n-\mu\cdot n}{\sigma\sqrt{n}}$ converges as a random 
variable to a normal (gaussian) law.
\end{theorem}
\begin{remark}
Again the theorem was originally proved for polynomial systems, but the arguments of the proof
hold more generally when $\Phi$ is rational. The role of the 
condition involving the existence of a nonsingular $3\times 3$
matrix is to grant the strict positivity of $\sigma$, as recently proved in~\cite{Ju10}. 
\end{remark}

\begin{table}
\begin{center}
\includegraphics[width=14cm]{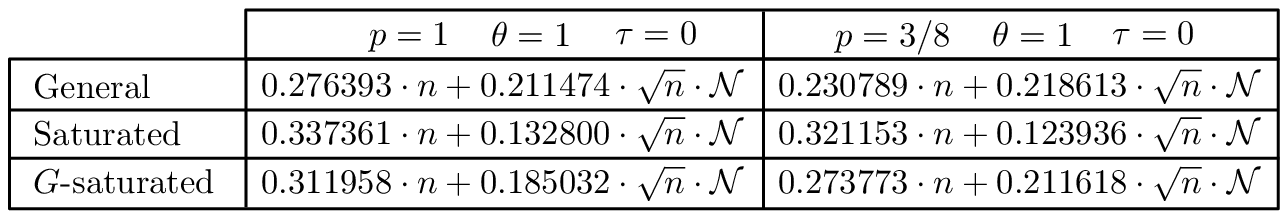}
\end{center}
\caption{Asymptotic behaviour of the number of links ($\mathcal{N}$ denotes a normal gaussian law).}
\label{fig:table_asympt2}
\end{table}

\begin{proposition}\label{prop:asympt2}
Let $p>0$. For $n\geq 1$, let $X_n$ be the number of links in a general (resp. saturated, resp. G-saturated) secondary structure  of length $n$
taken at random with weight proportional to $p^{\#(\mathrm{links})}$ (uniformly at random when $p=1$). 
Then there are computable strictly positive constants $\mu$ and $\sigma$ (depending on $p$, $\theta$, $\tau$, and on which 
setting: general, saturated, or G-saturated) such that 
$\frac{X_n-\mu\cdot n}{\sigma \sqrt{n}}$ converges as a random variable to a normal (gaussian) law.    
\end{proposition}
\proof
In each of the three settings (general, saturated, G-saturated), we have called $g(t,s)$ the enumerative generating function of
secondary structures with at least one link. We have seen that there are two rational expressions $Q(t,s,y)$ and $R(t,s,y)$ that  have  
nonnegative coefficients (in the series-expansion), and  there is an adjoint generating function $f(t,s)$ such that $f(t,s)=Q(t,s,f(t,s))$
and $g(t,s)=R(t,s,f(t,s))$; and the convergence domain of $Q(t,s,y)$
is the same as the convergence domain of $R(t,s,y)$. Note that the bivariate series $g(t,put^2)$ (with variables $t$ and $u$)
is the weighted generating function of secondary structures (with at least one link) where $t$ marks the length, $u$ marks the number of links,
and where each structure has weight $p^{\#(\mathrm{links})}$. 
It is easily checked that, if we set  
$a(t,u):=f(t,put^2)-\mathbf{1}_{\theta=0}$  (with $\theta$ the threshold parameter)  
 and $b(t):=g(t,put^2)$, 
then we are in the conditions of Theorem~\ref{theo:drmota}, 
with $\Phi(t,u,y):=Q(t,put^2,y+\mathbf{1}_{\theta=0})-\mathbf{1}_{\theta=0}$ and 
$\Psi(t,u,y):=R(t,put^2,y+\mathbf{1}_{\theta=0})$. Indeed the  $3\times 3$ matrix condition is readily checked, and  
for $u=1$ we get the equation of Proposition~\ref{prop:asympt1}, where we have already checked that the conditions are satisfied. 
\hfill $\Box$ \medskip

Table~\ref{fig:table_asympt2} shows the asymptotic behaviour for some standard parameter values. 
(The methodology to compute $\mu$ for saturated structures using computer algebra tools
is detailed in~\cite{Clote.jbcb09}.)   
The case $p=1$ corresponds to a homopolymer of length $n$ taken uniformly at random, while the case  $p=3/8$ corresponds to a (uniformly) 
random secondary structure where the underlying sequence is (any word) 
in $\{A,G,C,U\}^n$.  
As expected, saturated structures tend to have more links than G-saturated structures,
which tend to have more links than general structures.  

\section{Inclusion of dangles}
\label{section:eric3}

\begin{figure}
\begin{center}
\includegraphics[width=15cm]{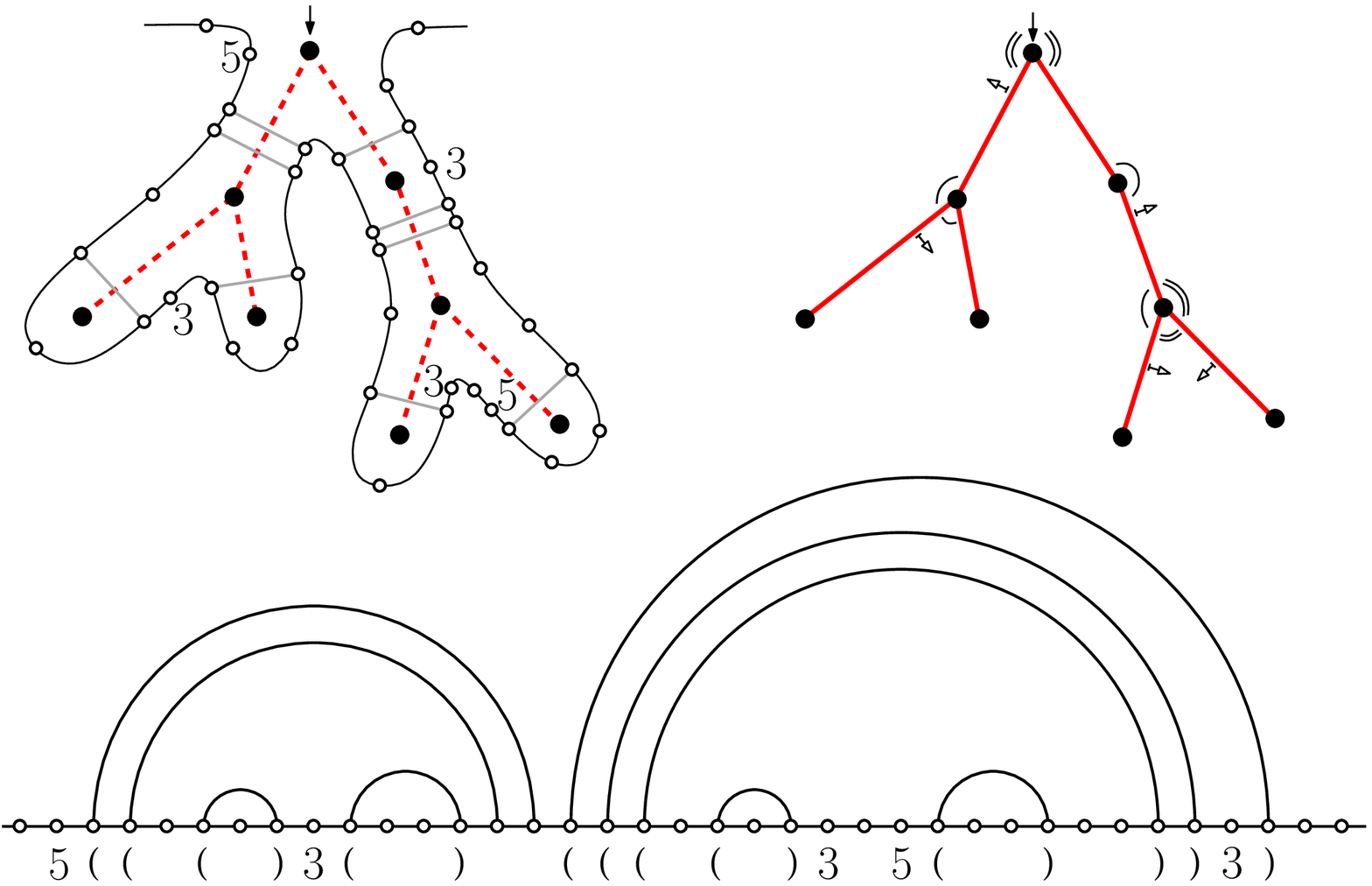}
\end{center}
\caption{Bottom: a secondary structure with dangles (two 
$5'$-dangles and three $3'$-dangles). Top-left:
the secondary structure with the dual plane tree. Top-right: each dangle yields a marked
edge-side in the dual plane tree (corners at inner nodes are simply marked if they have weight $1$,   are doubly marked if they have weight greater than $1$).}
\label{fig:dangles}
\end{figure}

We show here that the counting approach developed so far (based on duality with plane trees, generating functions, and substitution operations) can be easily adapted to take the presence
of so-called dangling bases into account.
In the parenthesis representation of the secondary structure (see Figure~\ref{fig:dangles}, bottom) 
a \emph{dangling base} (shortly dangle) is a distinguished free base of two possible kinds:
a \emph{$5$-dangle} has to be just before an opening parenthesis, a 
\emph{$3'$-dangle}
has to be just after of a closing parenthesis. Note that a dangling base that is just before
an opening parenthesis and just after a closing parenthesis is either a 
$5'$-dangle or a $3'$-dangle
but not both. For a structure with dangling bases, the \emph{underlying} secondary structure
is the structure where dangling bases are considered as usual free bases (i.e., are not
distinguished). 

\begin{figure}
\begin{center}
\includegraphics[width=15cm]{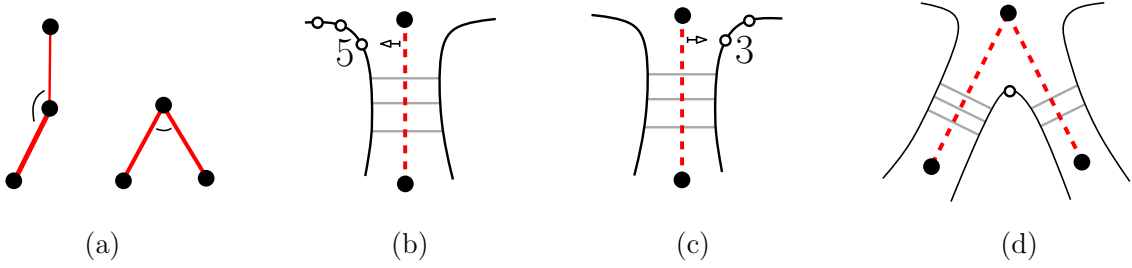}
\end{center}
\caption{(a) The first drawing shows a lateral corner, the second drawing shows an extremal corner 
(depending edges are bolder). (b) A $5'$-dangle yields a marked left-side of edge in the dual tree.
(c) A $3'$-dangle yields a marked right-side of edge in the dual tree. (d) 
Situation of  an extremal corner of weight 1, in which case 
 the two depending edge-sides can not both be marked.}
\label{fig:dangle_small}
\end{figure}

In the dual plane tree, a $5'$-dangle (resp. a $3'$-dangle) is indicated by a marked edge-side to the left
(resp. to the right) of the edge, see Figure~\ref{fig:dangle_small}(b)-(c).
To take dangles into account in our counting method, we need to distinguish two types
of corners in the dual plane tree $T$: a corner $c$ at vertex $v$ is called \emph{lateral} if
$c$ is incident to the edge going to the parent of $v$ in $T$ (when $v$ is not the root-vertex)
or $c$ is incident to the root (when $v$ is the root-node); note that every
inner node 
has two incident lateral corners (one on the left side, one on the right side).
The other corners at inner nodes in the tree are called \emph{extremal}, see Figure~\ref{fig:dangle_small}(a).
Given a corner $c$ of $T$ (at an inner node), 
an edge-side $s$ incident to $c$ is said to \emph{depend} on $c$ if 
$c$ is incident to $s$ at the extremity of $s$ closest to the root; 
note that a lateral corner has one depending edge-side while
an extremal corner has two depending edge-sides, see Figure~\ref{fig:dangle_small}(a).  

We now make important observations to determine when the 
 edge-sides depending on a corner $c$ can be marked. 
If $c$ has weight $0$ then none of the depending edge-sides can be marked,
because there is no free base in the sector of $c$ (hence no candidate to become a dangle). 
If $c$ is lateral and has positive weight (i.e., is a marked corner) then the depending
edge-side is allowed to be marked. If $c$ is extremal of weight $1$ then at most one
of the two edge-sides depending on $c$ is allowed to be marked (because the unique
free base in the sector of $c$ can not be both a $5'$-dangle and a $3'$-dangle).  
If $c$ is extremal of weight at least $2$ then the two depending edge-sides are allowed
to be marked (and are allowed to be both marked). 

Given these observations we can easily include a variable for dangles in our generating
function expressions (recall we have treated $3$ cases: general, saturated, G-saturated).  

\noindent{\bf General structures, inclusion of dangles in the results of Section~\ref{sec:GF}.} Denote by $F\equiv F(u,v_1,v_2,x)$ the generating
function of $\mathcal{F}$ (the one defined in Section~\ref{sec:GF}) where $u$ marks the number of leaves, $v_1$ (resp. $v_2$) marks
the number of marked corners that are lateral (resp. extremal), and $x$ marks the number of edges. 
Since a tree-vertex with $k\geq 1$ children has two incident corners that are lateral (the $k-1$ other
ones are extremal), we get the following equation (which specifies $F$ uniquely):
\begin{eqnarray*}
F&=&u+(2v_1+v_1^2)xF+\sum_{k\geq 2}x^k(1+v_1)^2(1+v_2)^{k-1}F^k\\
&=&u+\frac{x(1+v_1)^2F}{1-x(1+v_2)F}-xF. 
\end{eqnarray*}
Similarly, denoting by $G\equiv G(u,v_1,v_2,x)$ the generating function of $\mathcal{G}$ (where
the variables have the same meaning as for $F$), we have
$$
G=\frac{x(1+v_1)^2F}{1-x(1+v_2)F}.
$$
Let $g(t,s,r)$ be the generating function counting secondary structures with at least one link,
where $t$ marks the number of free elements (including dangles), $s$ marks the number of edges,
and $r$ marks the number of dangles. Then $g(t,s,r)=G(U,V_1,V_2,X)$, where
$$
U=\frac{t^{\theta}}{1-t},\ \ V_1=\frac{t(1+r)}{1-t},\ \ V_2=\frac{t(1+r)^2}{1-t}-tr^2,\ \ X=\frac{s^{\tau+1}}{1-s}.
$$
For $p>0, q\geq 0$, let $g_{p,q}(t)$ be the weighted generating function of secondary structures
where each structure has weight $p^{\#(\mathrm{links})}q^{\#(\mathrm{dangles})}$. 
Then $g_{p,q}(t)=g(t,pt^2,q)+t/(1-t)$. For instance, for $\theta=1$ and $\tau=0$ we find
$$g_{p,q}(t)=t+t^2+(1+p)t^3+(1+3p+2pq)t^4+(1+6p+p^2+6pq+pq^2)t^5+\cdots.$$

\vspace{.2cm}

\noindent{\bf Saturated structures, inclusion of dangles in the results of Section~\ref{sec:sat}.}
The equation for $F$ obtained in Section~\ref{sec:sat} becomes (when splitting $v$ into two variables
$v_1,v_2$ respectively for lateral and extremal marked corners):
\begin{eqnarray*}
F&=&u+2v_1xF+\sum_{k\geq 2}(1+2v_1+(k-1)v_2)x^kF^k\\
&=&u+\frac{x^2F^2+2xF\cdot(v_1-v_2)}{1-xF}+\frac{v_2}{(1-xF)^2}-v_2,
\end{eqnarray*}
and the expression of $G$ becomes
\begin{eqnarray*}
G&=&\sum_{k\geq 1}(1+2v_1+(k-1)v_2)x^kF^k\\
&=&\frac{xF\cdot\big(1+2(v_1-v_2)\big)}{1-xF}+\frac{v_2}{(1-xF)^2}-v_2.
\end{eqnarray*}
A structure with dangles is called \emph{saturated} if the underlying secondary structure
is saturated. 
Let $g(t,s,r)$ be the generating function counting saturated structures with at least one link,
where $t$ marks the number of free elements (including dangles), $s$ marks the number of edges,
and $r$ marks the number of dangles. Then $g(t,s,r)=G(U,V_1,V_2,X)$, where
$$
U=t^{\theta}(1+t),\ \ V_1=\frac{t-t^{\theta+2}}{1-t}(1+r),\ \ V_2=\frac{t-t^{\theta+2}}{1-t}(1+r)^2-tr^2,\ \ X=\frac{s}{1-s}.
$$
For $p>0, q\geq 0$, let $g_{p,q}(t)$ be the weighted generating function of saturated structures
where each structure has weight $p^{\#(\mathrm{links})}q^{\#(\mathrm{dangles})}$. 
Then $g_{p,q}(t)=g(t,pt^2,q)+t+\cdots+t^{\theta+1}$. For $\theta=1$ and $\tau=0$ we find
$$g_{p,q}(t)=t+t^2+pt^3+(3p+2pq)t^4+(4p+p^2+4pq)t^5+(2p+6p^2+2pq+4p^2q)t^6+\cdots.$$

\vspace{.2cm}

\noindent{\bf G-saturated structures, inclusion of dangles in the results of Section~\ref{sec:Gsat}.}
Let $C_k(v_1,v_2)$ be the polynomial generating function for
independent sets of the cycle $(1,\ldots,k)$, where $v_1$ (resp. $v_2$) marks the number of elements of the 
independent set that belong to $\{1,k\}$ (resp. to $\{2,\ldots,k-1\}$). Let 
$S_k(v)$ be the polynomial generating function for independent sets of the chain $1,\ldots,k$, where $v$ marks the number of elements in the independent set. Recall that $S(v,z):=\sum_{k\geq 0}S_k(v)z^k$ is given by
 $$S(v,z)=\frac{1+vz}{1-z-vz^2}.$$ 
Then one easily sees that for $k\geq 3$,
$$
C_k(v_1,v_2)=2v_1S_{k-3}(v_2)+S_{k-2}(v_2),
$$
and the equation for $F$ obtained in Section~\ref{sec:Gsat} becomes (when splitting $v$ into two variables
$v_1,v_2$ respectively for lateral and extremal marked corners):
\begin{eqnarray*}
F&=&u+2v_1xF+\sum_{k\geq 2}C_{k+1}(v_1,v_2)x^kF^k\\
&=&u+2v_1xF+\sum_{k\geq 2}(2v_1S_{k-2}(v_2)+S_{k-1}(v_2))x^kF^k\\
&=&u+2v_1xF+2v_1x^2F^2\cdot S(v_2,xF)+xF\cdot(S(v_2,xF)-1),
\end{eqnarray*}
which yields the simplified equation
$$
F=u+2v_1xF+\frac{1+2v_1x^2F^2\cdot(1+v_2xF)}{1-xF-v_2x^2F^2}-xF-1.
$$
And the expression of $G$ becomes at first
$$
G=\sum_{k\geq 1}\Big(S_{k-1}(v_2)+2v_1S_{k-2}(v_2)+v_1^2S_{k-3}(v_2)\Big)x^kF^k,
$$
with the conventions $S_{-1}(v)=1$, $S_{-2}(v)=0$. 
Hence we have
\begin{eqnarray*}
G&=&xF\cdot(1+v_1xF)^2S(v_2,xF)+2v_1xF+v_1^2x^2F^2\\
&=&\frac{xF\cdot \big(1+2v_1+xF\cdot(v_2+v_1^2)\big)}{1-xF-v_2x^2F^2}.
\end{eqnarray*}
A structure with dangles is called \emph{G-saturated} if the underlying secondary structure
is G-saturated. 
Let $g(t,s,r)$ be the generating function counting G-saturated structures with at least one link,
where $t$ marks the number of free elements (including dangles), $s$ marks the number of edges,
and $r$ marks the number of dangles. Then $g(t,s,r)=G(U,V_1,V_2,X)$, where
$$
U=t^{\theta}(1+t),\ \ V_1=\frac{t(1+r)}{1-t},\ \ V_2=\frac{t(1+r)^2}{1-t}-tr^2,\ \ X=\frac{s^{\tau+1}}{1-s}.
$$
For $p>0, q\geq 0$, let $g_{p,q}(t)$ be the weighted generating function of G-saturated structures
where each structure has weight $p^{\#(\mathrm{links})}q^{\#(\mathrm{dangles})}$. 
Then $g_{p,q}(t)=g(t,pt^2,q)+t/(1-t)$. For $\theta=1$ and $\tau=0$ we find
$$g_{p,q}=t+t^2+(1+p)t^3+(1+3p+2pq)t^4+(1+4p+p^2+4pq)t^5+(1+4p+6p^2+4pq+4p^2q)t^6+\cdots.$$

\vspace{.2cm}

\begin{table}
\begin{center}
\includegraphics[width=14cm]{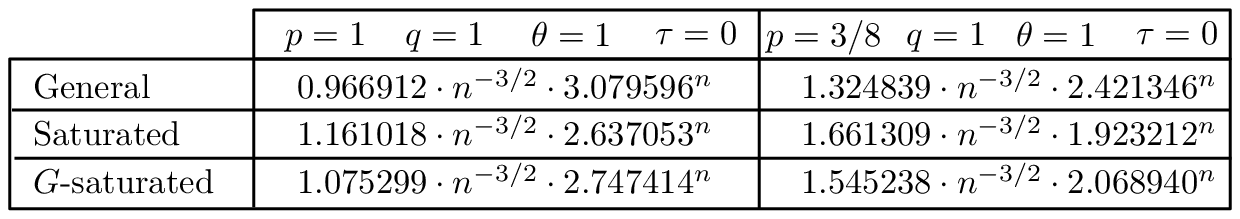}
\end{center}
\caption{Asymptotic behaviour of the $n$th coefficient  of the generating
function $g_{p,1}(t)$ counting secondary structures (general, saturated, or G-saturated) with dangles,  
with weight $p$ on each link.}
\label{fig:table_asympt_dangles}
\end{table}

\begin{table}
\begin{center}
\includegraphics[width=14cm]{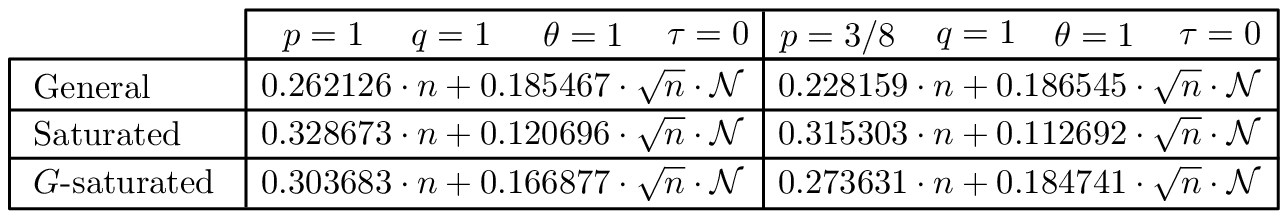}
\end{center}
\caption{Asymptotic behaviour of the number of links ($\mathcal{N}$ denotes a normal gaussian law)
for secondary structures (general, saturated, or G-saturated) with dangles, with weight $p$ 
on each link.}
\label{fig:table_asympt2_dangles}
\end{table}

\noindent{\bf Asymptotic results.} Propositions~\ref{prop:asympt1} and~\ref{prop:asympt2} directly extend to the case of any weight $q\geq 0$ for dangles (the case without dangles is $q=0$). 
 We give the numeric values corresponding to $q=1$ (asymptotic enumeration of structures
with dangles) in Tables~\ref{fig:table_asympt_dangles} and~\ref{fig:table_asympt2_dangles}, which are the counterparts of Tables~\ref{fig:table_asympt} and~\ref{fig:table_asympt2}.

\section{Discussion}
\label{section:discussion}

In this paper, we presented various context free grammars that generate
the set of secondary structures, according to different energy models:
Nussinov energy, base stacking energy, Turner energy,\footnote{Exact base
stacking parameters are ignored as is entropy; however, the context-free 
grammar allows the separate marking of distinct features, such as stacked
base pairs, hairpins, bulges, internal loops, multiloops.}
Turner with dangles (where dangles are rigorously
treated by the method of Markham and Zuker \cite{Markham.mmb08,markhamPhD}),
Turner (with external dangles), as well as saturated and G-saturated structures.
Using DSV, dominant singularity analysis and the Flajolet-Odlyzko theorem,
we proved that the asymptotic number
of secondary structures with annotated dangles, as computed in the partition
function of the Markham-Zuker software {\tt UNAFOLD} \cite{Markham.mmb08},
is $0.63998 \cdot n^{-3/2} \cdot  3.06039^n$, exponentially larger than the
number of all secondary structures
$1.104366 \cdot n^{-3/2} \cdot 2.618034^n$, previously established by
Stein and Waterman \cite{steinWaterman}. This result provides a partial
explanation for M.  Zuker's observation (personal communication) 
that {\tt UNAFOLD} requires substantially more computation time when dangles 
are included. 

Since the Nussinov energy model and the base stacking energy
model superficially appear to be almost equivalent, we 
presented a computational result that displays their marked 
differences.\footnote{Sheikh et al.  \cite{ponty:pseudoknotNPcomplete} 
show that minimum energy pseudoknotted structure prediction is NP-complete, in
contrast with the existence of a cubic time algorithm for the Nussinov
energy model \cite{tabaskaCaryGabowStormo}.} In particular,
the base stacking energy model leads to more cooperative folding and a
higher melting temperature for homopolymers than does the Nussinov energy
model. 

Finally, in the main part of the paper, we give generating
functions for the number of secondary structures and locally optimal
secondary structures, with respect to the Nussinov model and the 
base stacking energy models, permitting the determination of the
asymptotic number of (all resp. saturated resp. G-saturated) structures 
and the expected number of their base pairs, optionally requiring a
minimum stem length and stickiness parameter. With stickiness
parameter $2(p_{GC} + p_{AG} + p_{AU}) = \frac{3}{8}$, we 
obtain combinatorial results for RNA sequences using a reasonable
theoretical model. The principal advantage of our uniform treatment,
using duality, substitution of generating functions and the Drmota-Lalley-Woods
theorem is that with little additional effort, we can determine the
asymptotic number of (all resp. saturated resp. G-saturated) structures
with external dangles, and their expected number of base pairs. Such
computations would have been more difficult using grammars, DSV and
singularity analysis.

\section*{Acknowledgements}
Figure~\ref{fig:rnaSecStr2} was created by W.A. Lorenz and H. Jabbari.
We would like to thank the anonymous referees for their helpful comments.
\'E. Fusy is supported by the European project ExploreMaps ---ERC StG 208471. 
P. Clote is supported by the National Science Foundation under grants
DBI-0543506 and DMS-0817971, and by Digiteo Foundation project
{\em RNAomics}.  Any opinions, findings,
and conclusions or recommendations expressed in this material are
those of the authors and do not necessarily reflect the views of the
National Science Foundation.

\bibliographystyle{plain}

\end{document}